\documentclass{article}
\usepackage{amssymb}
\usepackage{amsfonts}
\usepackage{amsmath}

\input{tcilatex}
\begin{document}

\title{On some Gaussian Bernstein processes in $\mathbb{R}^{N}$ and the
periodic Ornstein-Uhlenbeck process}
\author{Pierre-A. Vuillermot$^{\ast }$ and Jean-C. Zambrini$^{\ast \ast }$ \\
UMR-CNRS 7502, Inst. \'{E}lie Cartan de Lorraine, Nancy, France$^{\ast }$\\
Departamento de Matem\'{a}tica, Universidade de Lisboa, Portugal$^{\ast \ast
}$}
\date{}
\maketitle

\begin{abstract}
In this article we prove new results regarding the existence of Bernstein
processes associated with the Cauchy problem of certain forward-backward
systems of decoupled linear deterministic parabolic equations defined in
Euclidean space of arbitrary dimension $N\in \mathbb{N}^{+}$, whose initial
and final conditions are positive measures. We concentrate primarily on the
case where the elliptic part of the parabolic operator is related to the
Hamiltonian of an isotropic system of quantum harmonic oscillators. In this
situation there are many Gaussian processes of interest whose existence
follows from our analysis, including $N$-dimensional stationary and
non-stationary Ornstein-Uhlenbeck processes, as well as a Bernstein bridge
which may be interpreted as a Markovian loop in a particular case. We also
introduce a new class of stationary non-Markovian processes which we
eventually relate to the $N$-dimensional periodic Ornstein-Uhlenbeck
process, and which is generated by a one-parameter family of non-Markovian
probability measures. In this case our construction requires the
consideration of an infinite hierarchy of pairs of forward-backward heat
equations associated with the pure point spectrum of the elliptic part,
rather than just one pair in the Markovian case. We finally stress the
potential relevance of these new processes to statistical mechanics, the
random evolution of loops and general pattern theory.
\end{abstract}

\section{Introduction and outline}

Let us consider the two adjoint parabolic Cauchy problems%
\begin{align}
\partial _{t}u(\mathsf{x},t)& =\frac{1}{2}\Delta _{\mathsf{x}}u(\mathsf{x}%
,t)-V(\mathsf{x})u(\mathsf{x},t),\text{ \ }(\mathsf{x},t)\in \mathbb{R}%
^{N}\times \left( 0,T\right] ,  \notag \\
u(\mathsf{x},0)& =\varphi _{0}(\mathsf{x}),\text{ \ \ }\mathsf{x}\in \mathbb{%
R}^{N}  \label{cauchyforward}
\end{align}%
and%
\begin{align}
-\partial _{t}v(\mathsf{x},t)& =\frac{1}{2}\Delta _{\mathsf{x}}v(\mathsf{x}%
,t)-V(\mathsf{x})v(\mathsf{x},t),\text{ \ }(\mathsf{x},t)\in \mathbb{R}%
^{N}\times \left[ 0,T\right) ,  \notag \\
v(\mathsf{x,}T)& =\psi _{T}(\mathsf{x}),\text{ \ \ }\mathsf{x}\in \mathbb{R}%
^{N},  \label{cauchybackward}
\end{align}%
where $T\in \left( 0,+\infty \right) $ is arbitrary. In these equations, $%
\Delta _{\mathsf{x}}$ denotes Laplace's operator with respect to the spatial
variable, $V$ is a real-valued function while $\varphi _{0}$ and $\psi _{T}$
are positive measures on $\mathbb{R}^{N}$. Both (\ref{cauchyforward}) and (%
\ref{cauchybackward}) can then be looked upon as defining a forward-backward
system of decoupled linear deterministic heat equations in Euclidean space.
To wit, the potential solutions to (\ref{cauchyforward}) wander off to the
future, whereas those of (\ref{cauchybackward}) evolve into the past. In
Section 2 below we show that under very general conditions on $V,\varphi
_{0} $ and $\psi _{T}$, we can associate with (\ref{cauchyforward})-(\ref%
{cauchybackward}) a class of the so-called Bernstein or reciprocal
processes, henceforth denoted by $Z_{\tau \in \left[ 0,T\right] }$. These
are processes that constitute a generalization of Markov processes which
have played an increasingly important r\^{o}le in various areas of
mathematics and mathematical physics over the years (see, e.g., \cite%
{albeyaza}, \cite{carlen}-\cite{cruzeirozambrini}, \cite{jamison}, \cite%
{roellythieullen}, \cite{vuizambrini} and the many references therein for a
history and ealier works on the subject). In the case of (\ref{cauchyforward}%
)-(\ref{cauchybackward}) the state space of the processes is the entire
Euclidean space, and their construction requires a transition function as
well as a joint probability distribution for $Z_{0}$ and $Z_{T}$, which we
denote respectively by $P$ and $\mu $ in the sequel. While $P$ depends
exclusively on the parabolic Green function associated with (\ref%
{cauchyforward})-(\ref{cauchybackward}), the measure $\mu $ involves both
Green's function and $\varphi _{0},\psi _{T}$, or more generally a
statistical mixture of $\varphi _{0}$ and $\psi _{T}$. It is $P$ and $\mu $
that allow one eventually to write out all the finite-dimensional
distributions of $Z_{\tau \in \left[ 0,T\right] }$, on which the remaining
part of this article is based. In Section 3 we concentrate primarily on
certain Markovian Gaussian Bernstein processes associated with a particular
case of (\ref{cauchyforward})-(\ref{cauchybackward}), namely,%
\begin{align}
\partial _{t}u(\mathsf{x},t)& =\frac{1}{2}\Delta _{\mathsf{x}}u(\mathsf{x}%
,t)-\frac{\lambda ^{2}}{2}\left\vert \mathsf{x}\right\vert ^{2}u(\mathsf{x}%
,t),\text{ \ }(\mathsf{x},t)\in \mathbb{R}^{N}\times \left( 0,T\right] , 
\notag \\
u(\mathsf{x},0)& =\varphi _{0}(\mathsf{x}),\text{ \ \ }\mathsf{x}\in \mathbb{%
R}^{N}  \label{forwardharmonic}
\end{align}%
and%
\begin{align}
-\partial _{t}v(\mathsf{x},t)& =\frac{1}{2}\Delta _{\mathsf{x}}v(\mathsf{x}%
,t)-\frac{\lambda ^{2}}{2}\left\vert \mathsf{x}\right\vert ^{2}v(\mathsf{x}%
,t),\text{ \ }(\mathsf{x},t)\in \mathbb{R}^{N}\times \left[ 0,T\right) , 
\notag \\
v(\mathsf{x,}T)& =\psi _{T}(\mathsf{x)},\text{ \ \ }\mathsf{x}\in \mathbb{R}%
^{N}  \label{backwardharmonic}
\end{align}%
where $\lambda >0$ and $\left\vert .\right\vert $ denotes the Euclidean norm
in $\mathbb{R}^{N}$, whose right-hand side is the Hamiltonian of an
isotropic system of quantum harmonic oscillators, up to a sign (see, e.g., 
\cite{messiah}). There we apply the results of Section 2 in order to
construct and analyze three processes of interest, by determining explicitly
in each case the finite-dimensional projections of the corresponding
Gaussian measures. The first one is stationary and thereby esssentially a $N$%
-dimensional Ornstein-Uhlenbeck process, while the second one is a
non-stationary Ornstein-Uhlenbeck process conditioned to start at the origin
of $\mathbb{R}^{N}$. The third one is a process that shares many of the
properties of a bridge, which we call a Bernstein bridge and which we may
identify with a Markovian loop in a particular case. The measures $\mu $ we
need for the construction of each one of these are intimately tied up with a
very specific choice of initial-final data in (\ref{forwardharmonic})-(\ref%
{backwardharmonic}). The situation is quite different in Section 4, where we
construct a new family of stationary non-Markovian Bernstein processes that
are related to the $N$-dimensional periodic Ornstein-Uhlenbeck process.
There, we prove that the relevant non-Markovian probability measures $\mu $
disintegrate into statistical mixtures of the form%
\begin{equation}
\mu =\sum_{\mathsf{m\in }\mathbb{N}^{N}}p_{\mathsf{m}}\mu _{\mathsf{m}}\text{%
, \ \ }p_{\mathsf{m}}>0\text{, \ \ }\sum_{\mathsf{m\in }\mathbb{N}^{N}}p_{%
\mathsf{m}}=1,  \label{statisticalmixture}
\end{equation}%
where each $\mu _{\mathsf{m}}$ is a measure related to initial-final
conditions $\varphi _{\mathsf{m,0}}$ and $\psi _{\mathsf{m,}T}$ in%
\begin{align}
\partial _{t}u(\mathsf{x},t)& =\frac{1}{2}\Delta _{\mathsf{x}}u(\mathsf{x}%
,t)-\frac{\lambda ^{2}}{2}\left\vert \mathsf{x}\right\vert ^{2}u(\mathsf{x}%
,t),\text{ \ }(\mathsf{x},t)\in \mathbb{R}^{N}\times \left( 0,T\right] , 
\notag \\
u(\mathsf{x},0)& =\varphi _{\mathsf{m,0}}(\mathsf{x}),\text{ \ \ }\mathsf{x}%
\in \mathbb{R}^{N}  \label{forwardharmonicbis}
\end{align}%
and%
\begin{align}
-\partial _{t}v(\mathsf{x},t)& =\frac{1}{2}\Delta _{\mathsf{x}}v(\mathsf{x}%
,t)-\frac{\lambda ^{2}}{2}\left\vert \mathsf{x}\right\vert ^{2}v(\mathsf{x}%
,t),\text{ \ }(\mathsf{x},t)\in \mathbb{R}^{N}\times \left[ 0,T\right) , 
\notag \\
v(\mathsf{x},T)& =\psi _{\mathsf{m,}T}(\mathsf{x}),\text{ \ \ }\mathsf{x}\in 
\mathbb{R}^{N},  \label{backwardharmonicbis}
\end{align}%
respectively. In other words, we show there that the construction of the $%
\mu $'s given by (\ref{statisticalmixture}) requires the consideration of a
hierarchy of infinitely many pairs of problems of the form (\ref%
{forwardharmonic})-(\ref{backwardharmonic}) associated with the whole pure
point spectrum of the elliptic operator on the right-hand side, rather than
just one pair in the Markovian case. Finally, we also point out the
potential applications of those new processes to statistical mechanics, to
the problem of random evolution of loops in space and to general pattern
theory, to name only three. In an appendix we also prove an important series
expansion for the Green function associated with(\ref{forwardharmonic})-(\ref%
{backwardharmonic}).

\section{A class of Bernstein processes in $\mathbb{R}^{N}$}

Generally speaking Bernstein processes can take values in any topological
space countable at infinity, and there are several equivalent ways to
characterize them (see, e.g., \cite{jamison}). However, the following
definition will be sufficient for our purposes:

\bigskip

\textbf{Definition 1.} Let $N\in \mathbb{N}^{+}$ and $T\in \left( 0,+\infty
\right) $ be arbitrary. We say the $\mathbb{R}^{N}$-valued process $Z_{\tau
\in \left[ 0,T\right] }$ defined on the complete probability space $\left(
\Omega ,\mathcal{F},\mathbb{P}\right) $ is a Bernstein process if 
\begin{equation}
\mathbb{E}\left( f(Z_{r})\left\vert \mathcal{F}_{s}^{+}\vee \mathcal{F}%
_{t}^{-}\right. \right) =\mathbb{E}\left( f(Z_{r})\left\vert
Z_{s},Z_{t}\right. \right)  \label{condiexpectations}
\end{equation}%
for every bounded Borel measurable function $f:\mathbb{R}^{N}\mapsto \mathbb{%
R}$, and for all $r,s,t$ satisfying $r\in \left( s,t\right) \subset \left[
0,T\right] $. In (\ref{condiexpectations}), the $\sigma $-algebras are%
\begin{equation}
\mathcal{F}_{s}^{+}=\sigma \left\{ Z_{\tau }^{-1}\left( E\right) :\tau \leq
s,\text{ }E\in \mathcal{B}_{N}\right\}  \label{pastalgebra}
\end{equation}%
and%
\begin{equation}
\mathcal{F}_{t}^{-}=\sigma \left\{ Z_{\tau }^{-1}\left( E\right) :\tau \geq
t,\text{ }E\in \mathcal{B}_{N}\right\} ,  \label{futurealgebra}
\end{equation}%
where $\mathcal{B}_{N}$ stands for the Borel $\sigma $-algebra on $\mathbb{R}%
^{N}$.

\bigskip

The dynamics of such a process at any time $r\in \left( s,t\right) $ are,
therefore, solely determined by the properties of the process at times $s$
and $t$, irrespective of its behavior prior to instant $s$ and after instant 
$t$. Of course, it is this fact that generalizes the Markov property.

In order to associate a mere Bernstein process with (\ref{cauchyforward})-(%
\ref{cauchybackward}) we now impose the following hypothesis, which regards
Green's function alone:

\bigskip

(H1) The measurable function $V:\mathbb{R}^{N}\mapsto \mathbb{R}$ is such
that the parabolic Green function $g$ associated with (\ref{cauchyforward})-(%
\ref{cauchybackward}) is jointly continuous in all variables and satisfies

\begin{equation}
g(\mathsf{x},t,\mathsf{y})>0  \label{positivity}
\end{equation}%
for all $\mathsf{x,y}\in \mathbb{R}^{N}$ and every $t\in \left( 0,T\right] $.

\bigskip

Let us now write $\mathcal{M}\left( \mathbb{R}^{N}\times \mathbb{R}^{N},%
\mathbb{C}\right) $ for the space of measures we are interested in, namely,
the topological dual of the Fr\'{e}chet space of all complex-valued,
compactly supported functions on $\mathbb{R}^{N}\times \mathbb{R}^{N}$
endowed with the usual locally convex topology (see, e.g., \cite{schwartz}).
Having (\ref{positivity}) at our disposal, let us introduce the functions%
\begin{equation}
p\left( \mathsf{x},t;\mathsf{z},r;\mathsf{y},s\right) :=\frac{g(\mathsf{x}%
,t-r,\mathsf{z})g(\mathsf{z},r-s,\mathsf{y})}{g(\mathsf{x},t-s,\mathsf{y})}
\label{bernsteindensity}
\end{equation}%
and%
\begin{equation}
P\left( \mathsf{x},t;E,r;\mathsf{y},s\right) :=\dint\limits_{E}d\mathsf{z}%
p\left( \mathsf{x},t;\mathsf{z},r;\mathsf{y},s\right)
\label{transitionfunction}
\end{equation}%
for every $E\in \mathcal{B}_{N}$, both being well defined and positive for
all $\mathsf{x},\mathsf{y},\mathsf{z}\in \mathbb{R}^{N}$ and all $r,s,t$
satisfying $r\in \left( s,t\right) \subset \left[ 0,T\right] $. For every $%
F\in \mathcal{B}_{N}\times \mathcal{B}_{N}$, let us also consider a positive
measure $\mu \in \mathcal{M}\left( \mathbb{R}^{N}\times \mathbb{R}^{N},%
\mathbb{C}\right) \mathbb{\ }$such that%
\begin{equation}
\mu \left( F\right) :=\int_{F}d\mu \left( \mathsf{x,y}\right)
\label{probabilitymeasure}
\end{equation}%
defines a probability measure on $\mathcal{B}_{N}\times \mathcal{B}_{N}$,
thus satisfying 
\begin{equation}
\int_{\mathbb{R}^{N}\times \mathbb{R}^{N}}d\mu \mathsf{\left( \mathsf{x,y}%
\right) =1.\label{normalization}}
\end{equation}%
The knowledge of both (\ref{transitionfunction}) and (\ref%
{probabilitymeasure}) then makes it possible to associate a Bernstein
process with (\ref{cauchyforward})-(\ref{cauchybackward}). The precise
statement is the following:

\bigskip

\textbf{Theorem 1.} \textit{Assume that Hypothesis (H1)\ holds, and let }$P$ 
\textit{and }$\mu $\textit{\ be given by (\ref{transitionfunction}) and (\ref%
{probabilitymeasure})-(\ref{normalization}), respectively. Then, there
exists a probability space }$\left( \Omega ,\mathcal{F},\mathbb{P}_{\mu
}\right) $ \textit{supporting an }$\mathbb{R}^{N}$\textit{-valued Bernstein
process }$Z_{\tau \in \left[ 0,T\right] }$\textit{\ such that the following
properties are valid:}

\textit{(a) The function }$P$\textit{\ is the transition function of }$%
Z_{\tau \in \left[ 0,T\right] }$\textit{\ in the sense that} 
\begin{equation*}
\mathbb{P}_{\mu }\left( Z_{r}\in E\left\vert Z_{s},Z_{t}\right. \right)
=P\left( Z_{t},t;E,r;Z_{s},s\right)
\end{equation*}%
\textit{for each }$E\in \mathcal{B}_{N}$ \textit{and all }$r,s,t$\textit{\
satisfying }$r\in \left( s,t\right) \subset \left[ 0,T\right] .$

\textit{(b) For every }$n\in \mathbb{N}^{+}$, \textit{the finite-dimensional
distributions of the process are given by}%
\begin{eqnarray}
&&\mathbb{P}_{\mu }\left( Z_{t_{1}}\in E_{1},...,Z_{t_{n}}\in E_{n}\right)  
\notag \\
&=&\int_{\mathbb{R}^{N}\times \mathbb{R}^{N}}\frac{d\mu \mathsf{\left( 
\mathsf{x,y}\right) }}{g(\mathsf{x},T,\mathsf{y})}\int_{E_{1}}d\mathsf{x}%
_{1}...\int_{E_{n}}d\mathsf{x}_{n}  \notag \\
&&\times \dprod\limits_{k=1}^{n}g\left( \mathsf{x}_{k},t_{k}-t_{k-1},\mathsf{%
x}_{k-1}\right) \times g\left( \mathsf{y},T-t_{n},\mathsf{x}_{n}\right) 
\label{distribution}
\end{eqnarray}%
\textit{for all }$E_{1},...,E_{n}\in \mathcal{B}_{N}$\textit{\ and all }$%
t_{0}=0<t_{1}<...<t_{n}<T$\textit{, where }$\mathsf{x}_{0}=\mathsf{x}$%
\textit{. In particular we have}%
\begin{eqnarray}
&&\mathbb{P}_{\mu }\left( Z_{t}\in E\right)   \notag \\
&=&\int_{\mathbb{R}^{N}\times \mathbb{R}^{N}}\frac{d\mu \mathsf{\left( 
\mathsf{x,y}\right) }}{g(\mathsf{x},T,\mathsf{y})}\int_{E}d\mathsf{z}g\left( 
\mathsf{x},t,\mathsf{z}\right) g\left( \mathsf{z},T-t,\mathsf{y}\right) 
\label{probability1}
\end{eqnarray}%
\textit{for each }$E\in \mathcal{B}_{N}$\textit{\ and every} $t\in \left(
0,T\right) $. \textit{Moreover,}%
\begin{equation}
\mathbb{P}_{\mu }\left( Z_{0}\in E\right) =\mu \left( E\times \mathbb{R}%
^{N}\right)   \label{probability2}
\end{equation}%
\textit{and}%
\begin{equation}
\mathbb{P}_{\mu }\left( Z_{T}\in E\right) =\mu \left( \mathbb{R}^{N}\times
E\right)   \label{probability3}
\end{equation}%
\textit{for each} $E\in \mathcal{B}_{N}$.

\textit{(c) }$\mathbb{P}_{\mu }$\textit{\ is the only probability measure
leading to the above properties.}

\bigskip

\textbf{Proof.} Up to minor technical details, a direct adaptation of the
method developed in Section 2 of \cite{vuizambrini} leads to%
\begin{eqnarray}
&&\mathbb{P}_{\mu }\left( Z_{0}\in E_{0},Z_{t_{1}}\in E_{1},...,Z_{t_{n}}\in
E_{n},Z_{T}\in E_{T}\right)   \notag \\
&=&\int_{E_{0}\times E_{T}}d\mu \mathsf{\left( \mathsf{x,y}\right) }%
\int_{E_{1}}d\mathsf{x}_{1}...\int_{E_{n}}d\mathsf{x}_{n}\dprod%
\limits_{k=1}^{n}p\left( \mathsf{y},T;\mathsf{x}_{k},t_{k};\mathsf{x}%
_{k-1},t_{k-1}\right)   \label{distributionbis}
\end{eqnarray}%
for all $E_{0},...,E_{T}\in B_{N}$\ and all $t_{0}=0<t_{1}<...<t_{n}<T$,
where $\mathsf{x}_{0}=\mathsf{x}$. In particular we have%
\begin{equation}
\mathbb{P}_{\mu }\left( Z_{0}\in E_{0},Z_{T}\in E_{T}\right) =\mu \left(
E_{0}\times E_{T}\right)   \label{jointdistribution}
\end{equation}%
for all $E_{0},E_{T}\in B_{N}$, that is, (\ref{probabilitymeasure}) is the
joint probability distribution of $Z_{0}$\ and $Z_{T}$. Now, from (\ref%
{bernsteindensity}) we obtain%
\begin{eqnarray*}
&&\dprod\limits_{k=1}^{n}p\left( \mathsf{y},T;\mathsf{x}_{k},t_{k};\mathsf{x}%
_{k-1},t_{k-1}\right)  \\
&=&\dprod\limits_{k=1}^{n}\frac{g(\mathsf{y},T-t_{k},\mathsf{x}_{k})g(%
\mathsf{x}_{k},t_{k}-t_{k-1},\mathsf{x}_{k-1})}{g(\mathsf{y},T-t_{k-1},%
\mathsf{x}_{k-1})} \\
&=&\frac{1}{g(\mathsf{x},T,\mathsf{y})}\dprod\limits_{k=1}^{n}g(\mathsf{x}%
_{k},t_{k}-t_{k-1},\mathsf{x}_{k-1})\times g(\mathsf{y},T-t_{n},\mathsf{x}%
_{n})
\end{eqnarray*}%
after $n-1$ cancellations of factors in order to obtain the second equality,
so that (\ref{distribution}) follows by choosing $E_{0}=E_{T}=$ $\mathbb{R}%
^{N}$ in (\ref{distributionbis}). Relations (\ref{probability2}) and (\ref%
{probability3}) are just a particular case of (\ref{jointdistribution}). \ \ 
$\blacksquare $ \ \ 

\bigskip

Of course, we can say more about $Z_{\tau \in \left[ 0,T\right] }$ if we
know more about $\mu $. First, $Z_{\tau \in \left[ 0,T\right] }$ is
Markovian if there exist positive measures $\nu _{0},\nu _{T}\in \mathcal{M}%
\left( \mathbb{R}^{N},\mathbb{C}\right) $ such that%
\begin{equation}
\mu \left( F\right) =\int_{F}d\left( \nu _{0}\otimes \nu _{T}\right) \left( 
\mathsf{x,y}\right) g(\mathsf{x},T,\mathsf{y})  \label{markovianmeasures}
\end{equation}%
for every $F\in \mathcal{B}_{N}\times \mathcal{B}_{N}$, with%
\begin{equation}
\int_{\mathbb{R}^{N}\times \mathbb{R}^{N}}d\left( \nu _{0}\otimes \nu
_{T}\right) \left( \mathsf{x,y}\right) g(\mathsf{x},T,\mathsf{y})=1
\label{normalizationbis}
\end{equation}%
in which case we also say that $\mu $ is Markovian. This result can be
traced back to the more general Theorem 3.1 in \cite{jamison}, and allows us
to make a closer connection between $Z_{\tau \in \left[ 0,T\right] }$ and (%
\ref{cauchyforward})-(\ref{cauchybackward}) provided we impose the following
hypothesis:

\bigskip

(H2) The measures $\varphi _{0},\psi _{T}\in \mathcal{M}\left( \mathbb{R}%
^{N},\mathbb{C}\right) $ in (\ref{cauchyforward})-(\ref{cauchybackward}) are
positive, and there exist a unique classical positive solution to (\ref%
{cauchyforward}) and a unique classical positive solution to (\ref%
{cauchybackward}), namely,%
\begin{equation}
u_{\varphi _{0}}(\mathsf{x},t)=\int_{\mathbb{R}^{N}}d\varphi _{0}(\mathsf{y)}%
g(\mathsf{x},t,\mathsf{y})  \label{forwardsolution}
\end{equation}%
and 
\begin{equation}
v_{\psi _{T}}(\mathsf{x},t)=\int_{\mathbb{R}^{N}}d\psi _{T}(\mathsf{y)}g(%
\mathsf{x},T-t,\mathsf{y})\mathsf{,\label{backwardsolution}}
\end{equation}%
respectively.

\bigskip

We then have the following consequence of Theorem 1:

\bigskip

\textbf{Corollary 1.}\textit{\ Assume that Hypotheses (H1)-(H2) hold, and
let us choose }$\mu $ \textit{of the form (\ref{markovianmeasures}) with }$%
\nu _{0}=\varphi _{0}$\textit{\ and }$\nu _{T}=\psi _{T}$\textit{. Then, the
Bernstein process }$Z_{\tau \in \left[ 0,T\right] }$ \textit{of Theorem 1 is
Markovian, and its finite-dimensional distributions are given by}%
\begin{eqnarray}
&&\mathbb{P}_{\mu }\left( Z_{t_{1}}\in E_{1},...,Z_{t_{n}}\in E_{n}\right) 
\label{distributionter} \\
&=&\int_{E_{1}}d\mathsf{x}_{1}...\int_{E_{n}}d\mathsf{x}_{n}\dprod%
\limits_{k=2}^{n}g\left( \mathsf{x}_{k},t_{k}-t_{k-1},\mathsf{x}%
_{k-1}\right) \times u_{\varphi _{0}}(\mathsf{x}_{1},t_{1})v_{\psi _{T}}(%
\mathsf{x}_{n},t_{n})  \notag
\end{eqnarray}%
\textit{for every integer }$n\geq 2$, \textit{all }$E_{1},...,E_{n}\in 
\mathcal{B}_{N}$\textit{\ and all }$t_{0}=0<t_{1}<...<t_{n}<T$\textit{,
where }$\mathsf{x}_{0}=\mathsf{x}$\textit{. Moreover we have}%
\begin{eqnarray}
&&\mathbb{P}_{\mu }\left( Z_{t}\in E\right)   \notag \\
&=&\int_{E}d\mathsf{x}u_{\varphi _{0}}\left( \mathsf{x},t\right) v_{\psi
_{T}}\left( \mathsf{x},t\right)   \label{probability4}
\end{eqnarray}%
\textit{for each }$E\in \mathcal{B}_{N}$\textit{\ and every} $t\in \left(
0,T\right) $. \textit{Finally, }%
\begin{equation}
\mathbb{P}_{\mu }\left( Z_{0}\in E\right) =\int_{E}d\varphi _{0}\left( 
\mathsf{x}\right) v_{\psi _{T}}\left( \mathsf{x},0\right) 
\label{probability5}
\end{equation}%
\textit{and}%
\begin{equation}
\mathbb{P}_{\mu }\left( Z_{T}\in E\right) =\int_{E}d\psi _{T}\left( \mathsf{x%
}\right) u_{\varphi _{0}}\left( \mathsf{x},T\right)   \label{probability6}
\end{equation}%
\textit{for each} $E\in \mathcal{B}_{N}$.

\bigskip

\textbf{Proof.} We first rewrite (\ref{distribution}) as%
\begin{eqnarray}
&&\mathbb{P}_{\mu }\left( Z_{t_{1}}\in E_{1},...,Z_{t_{n}}\in E_{n}\right) 
\notag \\
&=&\int_{\mathbb{R}^{N}\times \mathbb{R}^{N}}\frac{d\mu \mathsf{\left( 
\mathsf{x,y}\right) }}{g(\mathsf{x},T,\mathsf{y})}\int_{E_{1}}d\mathsf{x}%
_{1}...\int_{E_{n}}d\mathsf{x}_{n}  \notag \\
&&\times \dprod\limits_{k=2}^{n}g\left( \mathsf{x}_{k},t_{k}-t_{k-1},\mathsf{%
x}_{k-1}\right) \times g\left( \mathsf{x}_{1},t_{1},\mathsf{x}\right)
g\left( \mathsf{y},T-t_{n},\mathsf{x}_{n}\right)  \label{modifdistribution}
\end{eqnarray}%
and then substitute the choice of the measure $\mu $ into (\ref{probability3}%
) and (\ref{modifdistribution}), using (\ref{forwardsolution})-(\ref%
{backwardsolution}) along with the symmetry property of $g$ with respect to
the spatial variables. \ \ $\blacksquare $

\bigskip

\textsc{Remarks.} (1) In the Markovian case, we have thus exhibited a very
general class of initial-final conditions in (\ref{cauchyforward})-(\ref%
{cauchybackward}) which allows us to determine the processes $Z_{\tau \in %
\left[ 0,T\right] }$ completely, including their marginal distributions (\ref%
{probability5})-(\ref{probability6}). The converse point of view was
developed in \cite{beurling} and its references, where it was shown instead
that it is a general class of marginal distributions which determines the
initial-final data of the relevant partial differential equations, through a
system of nonlinear integral equations. However, the resulting initial-final
conditions of \cite{beurling} belong to the class of positive continuous
functions, and not to the larger class of positive measures as is the case
in this article.

(2) A glance at (\ref{backwardsolution}) shows that it is sufficient to use
elementary time reversal in Green's function to obtain the solution to (\ref%
{cauchybackward}). Although the situation is not always that simple,
particularly when the given parabolic equations are non-autonomous, it is
still possible to define a quite appropriate probabilistic notion of time
symmetry in general. We refer the reader to \cite{vuizambrini} for further
details.

\bigskip

Relations (\ref{distribution}) and (\ref{distributionter}) are the
fundamental relations that will allow us to construct the Gaussian processes
associated with (\ref{forwardharmonic})-(\ref{backwardharmonic}) in the next
sections.

\section{Two Ornstein-Uhlenbeck processes and a Bernstein bridge}

In the remaining part of this article we denote by $\left( .,.\right) _{%
\mathbb{R}^{N}}$ the Euclidean inner product in $\mathbb{R}^{N}$, and by $%
L^{2}\left( \mathbb{R}^{N},\mathbb{C}\right) $ the usual Lebesgue space of
all complex-valued, square-integrable functions on $\mathbb{R}^{N}$. We
begin by considering the forward-backward system (\ref{forwardharmonic})-(%
\ref{backwardharmonic}) with centered Gaussian initial-final data, namely,%
\begin{eqnarray}
\partial _{t}u(\mathsf{x},t) &=&\frac{1}{2}\Delta _{\mathsf{x}}u(\mathsf{x}%
,t)-\frac{\lambda ^{2}\left\vert \mathsf{x}\right\vert ^{2}}{2}u(\mathsf{x}%
,t),\text{ \ \ }(\mathsf{x},t)\in \mathbb{R}^{N}\mathbb{\times }\left( 0,T%
\right] ,  \notag \\
u(\mathsf{x},0) &=&\varphi _{0,\lambda }(\mathsf{x})=\left( \frac{\lambda
\exp \left[ \lambda T\right] }{\pi }\right) ^{\frac{N}{4}}\exp \left[ -\frac{%
\lambda \left\vert \mathsf{x}\right\vert ^{2}}{2}\right]  \label{initial2}
\end{eqnarray}%
and%
\begin{eqnarray}
-\partial _{t}v(\mathsf{x},t) &=&\frac{1}{2}\Delta _{\mathsf{x}}v(\mathsf{x}%
,t)-\frac{\lambda ^{2}\left\vert \mathsf{x}\right\vert ^{2}}{2}v(\mathsf{x}%
,t),\text{ \ \ }(\mathsf{x},t)\in \mathbb{R}^{N}\mathbb{\times }\left[
0,T\right) ,  \notag \\
v(\mathsf{x},T) &=&\psi _{T,\lambda }(\mathsf{x})=\left( \frac{\lambda \exp %
\left[ \lambda T\right] }{\pi }\right) ^{\frac{N}{4}}\exp \left[ -\frac{%
\lambda \left\vert \mathsf{x}\right\vert ^{2}}{2}\right] ,  \label{final2}
\end{eqnarray}%
thereby identifying the measure $\varphi _{0,\lambda }=$ $\psi _{T,\lambda }$
with its Gaussian density relative to the Lebesgue measure in $\mathbb{R}%
^{N} $. Let us recall that the self-adjoint realization in $L^{2}\left( 
\mathbb{R}^{N},\mathbb{C}\right) $ of the elliptic operator on the
right-hand side of these equations has a pure point spectrum. More
specifically, for every $m\in \mathbb{N}$ let%
\begin{equation}
h_{m,\lambda }(x):=\sqrt[4]{\lambda }h_{m}\left( \sqrt{\lambda }x\right)
\label{scaledhermitefunctions}
\end{equation}%
be the one-dimensional, suitably scaled Hermite functions where%
\begin{equation}
h_{m}\left( x\right) =\left( \pi ^{\frac{1}{2}}2^{m}m!\right) ^{-\frac{1}{2}%
}e^{-\frac{x^{2}}{2}}H_{m}(x),  \label{hermitefunctions}
\end{equation}%
and where the $H_{m}$'s are the Hermite polynomials%
\begin{equation}
H_{m}(x)=(-1)^{m}e^{x^{2}}\frac{d^{m}}{dx^{m}}e^{-x^{2}}.
\label{hermitepolynomials}
\end{equation}%
The following spectral result regarding the operator on the right-hand side
of (\ref{initial2}) is well-known and in any case can be verified directly
by means of an explicit computation:

\bigskip

\textbf{Lemma 1.} \textit{The tensor products }$\otimes
_{j=1}^{N}h_{m_{j},\lambda }$\textit{\ where the }$m_{j}$'s\textit{\ vary
independently on }$\mathbb{N}$\textit{\ provide an orthonormal basis of }$%
L^{2}\left( \mathbb{R}^{N},\mathbb{C}\right) $\textit{, and moreover the
eigenvalue equation}%
\begin{equation}
\left( -\frac{1}{2}\Delta _{\mathsf{x}}+\frac{\lambda ^{2}\left\vert \mathsf{%
x}\right\vert ^{2}}{2}\right) \otimes _{j=1}^{N}h_{m_{j},\lambda }\left( 
\mathsf{x}\right) =\lambda \left( \sum_{j=1}^{N}m_{j}+\frac{N}{2}\right)
\otimes _{j=1}^{N}h_{m_{j},\lambda }\left( \mathsf{x}\right)
\label{eigenequation}
\end{equation}%
\textit{holds for every }$\mathsf{x\in }\mathbb{R}^{N}$.

\bigskip

In particular, by reference to (\ref{hermitefunctions}) and (\ref%
{hermitepolynomials}) we have%
\begin{equation}
\otimes _{j=1}^{N}h_{0,\lambda }\left( \mathsf{x}\right) =\left( \frac{%
\lambda }{\pi }\right) ^{\frac{N}{4}}\exp \left[ -\frac{\lambda \left\vert 
\mathsf{x}\right\vert ^{2}}{2}\right] >0  \label{groundstate}
\end{equation}%
whose associated eigenvalue is $E_{\lambda }=\frac{N\lambda }{2}$, so that (%
\ref{groundstate}) corresponds to the initial-final conditions in (\ref%
{initial2})-(\ref{final2}) up to a normalization factor chosen in such a way
that (\ref{normalization}) holds for the measure%
\begin{equation*}
\mu _{\lambda }\left( F\right) =\int_{F}d\mathsf{x}d\mathsf{y}\varphi
_{0,\lambda }\left( \mathsf{x}\right) \psi _{T,\lambda }\left( \mathsf{y}%
\right) g_{\lambda }(\mathsf{x},T,\mathsf{y})
\end{equation*}%
where $F\in \mathcal{B}_{N}\times \mathcal{B}_{N}$, and where 
\begin{eqnarray}
&&g_{\lambda }(\mathsf{x},t,\mathsf{y})  \label{Mehlerkernel} \\
&=&\left( \frac{\lambda }{2\pi \sinh \left( \lambda t\right) }\right) ^{%
\frac{N}{2}}\exp \left[ -\frac{\lambda \left( \cosh (\lambda t)\left(
\left\vert \mathsf{x}\right\vert ^{2}+\left\vert \mathsf{y}\right\vert
^{2}\right) -2\left( \mathsf{x,y}\right) _{\mathbb{R}^{N}}\right) }{2\sinh
\left( \lambda t\right) }\right]  \notag
\end{eqnarray}%
is the $N$-dimensional version of Mehler's kernel for $t\in \left( 0,T\right]
$ (see the appendix for details). In this way, the unique classical positive
solutions to (\ref{initial2})-(\ref{final2}) satisfying the requirement%
\begin{equation}
\int_{\mathbb{R}^{N}}d\mathsf{x}u_{\lambda }(\mathsf{x},t)v_{\lambda }(%
\mathsf{x},t)=1  \label{requirement}
\end{equation}%
for every $t\in \left[ 0,T\right] $ are%
\begin{equation}
u_{\lambda }(\mathsf{x},t)=\left( \frac{\lambda \exp \left[ \lambda T\right] 
}{\pi }\right) ^{\frac{N}{4}}\exp \left[ -\frac{\lambda \left( \left\vert 
\mathsf{x}\right\vert ^{2}+Nt\right) }{2}\right]  \label{forwardsolutionbis}
\end{equation}%
and 
\begin{equation}
v_{\lambda }(\mathsf{x},t)=\left( \frac{\lambda \exp \left[ \lambda T\right] 
}{\pi }\right) ^{\frac{N}{4}}\exp \left[ -\frac{\lambda \left( \left\vert 
\mathsf{x}\right\vert ^{2}+N\left( T-t\right) \right) }{2}\right] ,
\label{backwardsolutionbis}
\end{equation}%
respectively. Then the following result holds, where $\mathbb{E}_{\mu }$
stands for the expectation functional on $\left( \Omega ,\mathcal{F},\mathbb{%
P}_{\mu }\right) $:

\bigskip

\textbf{Proposition 1. }\textit{The Bernstein process }$Z_{\tau \in \left[
0,T\right] }^{\lambda }$\textit{\ associated with (\ref{initial2})-(\ref%
{final2}) in the sense of Corollary 1 is a Gaussian and Markovian process
such that}%
\begin{equation}
\mathbb{P}_{\mu }\left( Z_{t}^{\lambda }\in E\right) =(2\pi \sigma _{\lambda
})^{-\frac{N}{2}}\int_{E}d\mathsf{x}\exp \left[ -\frac{\left\vert \mathsf{x}%
\right\vert ^{2}}{2\sigma _{\lambda }}\right]  \label{gaussianprocess}
\end{equation}%
\textit{for each }$t\in \left[ 0,T\right] $ \textit{and every} $E\in 
\mathcal{B}_{N}$, \textit{where}%
\begin{equation*}
\sigma _{\lambda }=\frac{1}{2\lambda }.
\end{equation*}%
\textit{Furthermore, the components of }$Z_{\tau \in \left[ 0,T\right]
}^{\lambda }$\textit{\ satisfy the relation}%
\begin{equation}
\mathbb{E}_{\mu }\left( Z_{s}^{\lambda ,i}Z_{t}^{\lambda ,j}\right) =\frac{%
\exp \left[ -\lambda \left\vert t-s\right\vert \right] }{2\lambda }\delta
_{i,j}  \label{ornsteinuhlenbeck}
\end{equation}%
\textit{for all }$s,t\in \left[ 0,T\right] $\textit{\ and all }$i,j\in
\left\{ 1,...,N\right\} $\textit{. Thus, }$Z_{\tau \in \left[ 0,T\right]
}^{\lambda }$\textit{\ identifies in law with a process} \textit{whose
components are all independent, one-dimensional and stationary
Ornstein-Uhlenbeck processes.}

\bigskip

\textbf{Proof.} The fact that $Z_{\tau \in \left[ 0,T\right] }^{\lambda }$
is Gaussian and Markovian satisfying (\ref{gaussianprocess}) follows from
Corollary 1 with (\ref{Mehlerkernel}) and (\ref{forwardsolutionbis})-(\ref%
{backwardsolutionbis}) plugged into (\ref{distributionter}) and (\ref%
{probability4}), so that the density of the probability distribution for $%
\left( Z_{t_{1}}^{\lambda },...,Z_{t_{n}}^{\lambda }\right) \in \mathbb{R}%
^{nN}$ is%
\begin{eqnarray*}
&&\dprod\limits_{k=2}^{n}g_{\lambda }\left( \mathsf{x}_{k},t_{k}-t_{k-1},%
\mathsf{x}_{k-1}\right) \times u_{\lambda }(\mathsf{x}_{1},t_{1})v_{\lambda
}(\mathsf{x}_{n},t_{n}) \\
&=&\left( 2\pi \right) ^{-\frac{nN}{2}}\left( 2\lambda ^{n}e^{\lambda
(t_{n}-t_{1})}\right) ^{\frac{N}{2}}\left( \dprod\limits_{k=2}^{n}\sinh
\left( \lambda (t_{k}-t_{k-1})\right) \right) ^{-\frac{N}{2}} \\
&&\times \exp \left[ -\frac{\lambda }{2}\sum_{k=2}^{n}\frac{\cosh (\lambda
(t_{k}-t_{k-1}))\left( \left\vert \mathsf{x}_{k}\right\vert ^{2}+\left\vert 
\mathsf{x}_{k-1}\right\vert ^{2}\right) -2\left( \mathsf{x}_{k}\mathsf{,x}%
_{k-1}\right) _{\mathbb{R}^{N}}}{\sinh \left( \lambda (t_{k}-t_{k-1})\right) 
}\right] \\
&&\times \exp \left[ -\frac{\lambda }{2}\left( \left\vert \mathsf{x}%
_{1}\right\vert ^{2}+\left\vert \mathsf{x}_{n}\right\vert ^{2}\right) \right]
.
\end{eqnarray*}%
The $nN\times nN$ corresponding covariance matrix is then of the form $%
C_{\lambda }\otimes \mathbb{I}_{N}$ with $\mathbb{I}_{N}$ the identity
matrix in $\mathbb{R}^{N}$, where $C_{\lambda }^{-1}$ is tridiagonal and
obtained by identification of the quadratic form in the argument of the
above exponentials when $N=1$. This gives 
\begin{equation*}
C_{\lambda ,k,k}^{-1}=\left\{ 
\begin{array}{c}
\frac{\lambda e^{\lambda \left( t_{2}-t_{1}\right) }}{\sinh \left( \lambda
(t_{2}-t_{1})\right) }\text{ \ \ \ \ for }k=1, \\ 
\\ 
\frac{\lambda \sinh \left( \lambda (t_{k+1}-t_{k-1})\right) }{\sinh \left(
\lambda (t_{k+1}-t_{k})\right) \sinh \left( \lambda (t_{k}-t_{k-1})\right) }%
\text{ \ \ \ \ for }k=2,...,n-1, \\ 
\\ 
\frac{\lambda e^{\lambda \left( t_{n}-t_{n-1}\right) }}{\sinh \left( \lambda
(t_{n}-t_{n-1})\right) }\text{ \ \ \ \ for }k=n\text{ }%
\end{array}%
\right.
\end{equation*}%
(the second line not being there when $n=2)$, and%
\begin{equation*}
C_{\lambda ,k,k-1}^{-1}=C_{\lambda ,k-1,k}^{-1}=-\frac{\lambda }{\sinh
\left( \lambda \left\vert t_{k}-t_{k-1}\right\vert \right) }\text{ \ \ for }%
k=2,...,n.
\end{equation*}%
Consequently, using for instance the analytical inversion formulae in \cite%
{huangmcoll}, or by direct verification, we obtain%
\begin{equation*}
C_{\lambda ,k,l}=\frac{\exp \left[ -\lambda \left\vert
t_{k}-t_{l}\right\vert \right] }{2\lambda }
\end{equation*}%
for all $k,l\in \left\{ 1,...,n\right\} ,$ so that (\ref{ornsteinuhlenbeck})
eventually holds. Now, let us consider the forward It\^{o} integral equation 
\begin{eqnarray}
X_{t} &=&e^{-\lambda t}Z_{0}^{\lambda }+\int_{0}^{t}e^{-\lambda (t-\tau
)}dW_{\tau },\text{ \ \ }t\in \left[ 0,T\right] ,  \notag \\
X_{0} &=&Z_{0}^{\lambda }  \label{itoequation}
\end{eqnarray}%
where $W_{\tau \in \left[ 0,T\right] }$ is a given Wiener process in $%
\mathbb{R}^{N}$, and where $Z_{0}^{\lambda }$ is distributed according to (%
\ref{gaussianprocess}) and independent of $W_{\tau \in \left[ 0,T\right] }$.
It is well known that the solution process $X_{\tau \in \left[ 0,T\right]
}^{\lambda }$ to (\ref{itoequation}) is centered Gaussian with covariance (%
\ref{ornsteinuhlenbeck}) (see, e.g., Section 5.6 in Chapter 5 of \cite%
{karashreve}), so that $Z_{\tau \in \left[ 0,T\right] }^{\lambda }$
identifies in law with $X_{\tau \in \left[ 0,T\right] }^{\lambda }$.
Therefore, $Z_{\tau \in \left[ 0,T\right] }^{\lambda }$ is indeed a $N$%
-dimensional Ornstein-Uhlenbeck process with the stated properties. \ \ $%
\blacksquare $

\bigskip

Next, we show that if we require instead $Z_{0}$ to be a given point in $%
\mathbb{R}^{N}$, the solution process to (\ref{itoequation}) is also a
particular Bernstein process. We do this in the simplest case where $Z_{0}$
is the origin, by considering the forward-backward system%
\begin{align}
\partial _{t}u(\mathsf{x},t)& =\frac{1}{2}\Delta _{\mathsf{x}}u(\mathsf{x}%
,t)-\frac{\lambda ^{2}\left\vert \mathsf{x}\right\vert ^{2}}{2}u(\mathsf{x}%
,t),\text{ \ \ }(\mathsf{x},t)\in \mathbb{R}^{N}\mathbb{\times }\left( 0,T%
\right] ,  \notag \\
u(\mathsf{x},0)& =\varphi _{0,\lambda }(\mathsf{x})=\left( \exp \left[
\lambda T\right] \right) ^{\frac{N}{4}}\delta (\mathsf{x})
\label{initialvalueproblembis}
\end{align}%
and%
\begin{align}
-\partial _{t}v(\mathsf{x},t)& =\frac{1}{2}\Delta _{\mathsf{x}}v(\mathsf{x}%
,t)-\frac{\lambda ^{2}\left\vert \mathsf{x}\right\vert ^{2}}{2}v(\mathsf{x}%
,t),\text{ \ \ }(\mathsf{x},t)\in \mathbb{R}^{N}\mathbb{\times }\left[
0,T\right) ,  \notag \\
v(\mathsf{x},T)& =\psi _{T,\lambda }(\mathsf{x})=\left( \exp \left[ \lambda T%
\right] \right) ^{\frac{N}{4}}\exp \left[ -\frac{\lambda \left\vert \mathsf{x%
}\right\vert ^{2}}{2}\right] ,  \label{finalvalueproblembis}
\end{align}%
with $\delta $ the Dirac measure. The corresponding Markovian measure (\ref%
{probabilitymeasure}) is then determined by%
\begin{equation*}
\mu _{\lambda }\mathsf{\left( \mathsf{x,y}\right) =}\left( \exp \left[
\lambda T\right] \right) ^{\frac{N}{2}}\delta (\mathsf{x})\exp \left[ -\frac{%
\lambda \left\vert \mathsf{y}\right\vert ^{2}}{2}\right] g(\mathsf{x},T,%
\mathsf{y})
\end{equation*}%
and satisfies (\ref{normalization}), while the unique relevant classical
positive solutions to (\ref{initialvalueproblembis})-(\ref%
{finalvalueproblembis}) are%
\begin{equation}
u_{\lambda }(\mathsf{x},t)=\left( \frac{\lambda \exp \left[ \frac{\lambda T}{%
2}\right] }{2\pi \sinh (\lambda t)}\right) ^{\frac{N}{2}}\exp \left[ -\frac{%
\lambda \coth (\lambda t)\left\vert \mathsf{x}\right\vert ^{2}}{2}\right]
\label{forwardsolutionquarto}
\end{equation}%
and%
\begin{equation}
v_{\lambda }(\mathsf{x},t)=\left( \exp \left[ \lambda T\right] \right) ^{%
\frac{N}{4}}\exp \left[ -\frac{\lambda \left( \left\vert \mathsf{x}%
\right\vert ^{2}+N\left( T-t\right) \right) }{2}\right] ,
\label{finalsolutionquarto}
\end{equation}%
respectively. We then have the following result:

\bigskip

\textbf{Proposition 2. }\textit{The Bernstein process }$Z_{\tau \in \left[
0,T\right] }^{\lambda }$\textit{\ associated with (\ref%
{initialvalueproblembis})-(\ref{finalvalueproblembis}) in the sense of
Corollary 1 is a Gaussian and Markovian process such that}%
\begin{equation}
\mathbb{P}_{\mu }\left( Z_{t}^{\lambda }\in E\right) =(2\pi \sigma _{\lambda
}(t))^{-\frac{N}{2}}\int_{E}d\mathsf{x}\exp \left[ -\frac{\left\vert \mathsf{%
x}\right\vert ^{2}}{2\sigma _{\lambda }(t)}\right]
\label{gaussianprocessquarto}
\end{equation}%
\textit{for each }$t\in \left( 0,T\right] $\textit{\ and every} $E\in 
\mathcal{B}_{N}$, \textit{where}%
\begin{equation}
\sigma _{\lambda }(t)=\frac{\sinh (\lambda t)\exp \left[ -\lambda t\right] }{%
\lambda }.  \label{varianceter}
\end{equation}%
\textit{Furthermore we have}%
\begin{equation}
\mathbb{P}_{\mu }\left( Z_{0}^{\lambda }=\mathsf{o}\right) =1,
\label{probability8}
\end{equation}%
\textit{and} \textit{the components of }$Z_{\tau \in \left[ 0,T\right]
}^{\lambda }$\textit{\ satisfy the relation }%
\begin{equation}
\mathbb{E}_{\mu }\left( Z_{s}^{\lambda ,i}Z_{t}^{\lambda ,j}\right) =\frac{%
\exp \left[ -\lambda \left( t+s\right) \right] }{2\lambda }\left( \exp \left[
2\lambda \left( t\wedge s\right) \right] -1\right) \delta _{i,j}
\label{covarianceter}
\end{equation}%
\textit{for all }$s,t\in \left[ 0,T\right] $\textit{\ and all }$i,j\in
\left\{ 1,...,N\right\} $\textit{. Thus, }$Z_{\tau \in \left[ 0,T\right]
}^{\lambda }$\textit{\ identifies in law with a process whose components are
all independent, one-dimensional non-stationary Ornstein-Uhlenbeck processes
conditioned to start at the origin.}

\bigskip

\textbf{Proof.} As in the proof of Proposition 1, the first part of the
statement including (\ref{gaussianprocessquarto})-(\ref{probability8})
follows from the appropriate substitutions into the formulae of Corollary 1.
Furthermore, the matrix $C_{\lambda }^{-1}$ resulting from the
identification of the quadratic form in the Gaussian density of $\left(
Z_{t_{1}}^{\lambda },...,Z_{t_{n}}^{\lambda }\right) $ is the same as that
in Proposition 1, with the exception of $C_{\lambda 1,1}^{-1}$ which now
reads%
\begin{equation*}
C_{\lambda ,1,1}^{-1}=\frac{\lambda \sinh \left( \lambda t_{2}\right) }{%
\sinh \left( \lambda (t_{2}-t_{1})\right) \sinh \left( \lambda t_{1}\right) }%
\text{ }.\text{\ }
\end{equation*}%
Inverting again we eventually obtain%
\begin{equation*}
C_{\lambda ,k,l}=\frac{\exp \left[ -\lambda \left( t_{k}+t_{l}\right) \right]
}{2\lambda }\left( \exp \left[ 2\lambda \left( t_{k}\wedge t_{l}\right) %
\right] -1\right)
\end{equation*}%
for all $k,l\in \left\{ 1,...,n\right\} ,$ so that (\ref{covarianceter})
holds. Therefore, $Z_{\tau \in \left[ 0,T\right] }^{\lambda }$ is indeed a $%
N $-dimensional Ornstein-Uhlenbeck process with the stated properties. \ \ $%
\blacksquare $

\bigskip

\textsc{Remarks.} (1) It is equally easy to condition the Ornstein-Uhlenbeck
process so that it ends at the origin of $\mathbb{R}^{N}$ when $t=T$. The
underlying Bernstein process $\hat{Z}_{\tau \in \left[ 0,T\right] }^{\lambda
}$ is then simply determined by swaping the initial-final conditions in (\ref%
{initialvalueproblembis}) and (\ref{finalvalueproblembis}). Indeed, in doing
so the relevant solutions become%
\begin{equation*}
u_{\lambda }(\mathsf{x},t)=\left( \exp \left[ \lambda T\right] \right) ^{%
\frac{N}{4}}\exp \left[ -\frac{\lambda \left( \left\vert \mathsf{x}%
\right\vert ^{2}+Nt\right) }{2}\right]
\end{equation*}%
and%
\begin{equation*}
v_{\lambda }(\mathsf{x},t)=\left( \frac{\lambda \exp \left[ \frac{\lambda T}{%
2}\right] }{2\pi \sinh (\lambda (T-t))}\right) ^{\frac{N}{2}}\exp \left[ -%
\frac{\lambda \coth (\lambda (T-t))\left\vert \mathsf{x}\right\vert ^{2}}{2}%
\right]
\end{equation*}%
instead of (\ref{forwardsolutionquarto}) and (\ref{finalsolutionquarto}),
respectively, so that $\hat{Z}_{\tau \in \left[ 0,T\right] }^{\lambda }$ is
just the time-reversal of the process of Proposition 2, namely,%
\begin{equation*}
\hat{Z}_{\tau }^{\lambda }=Z_{T-\tau }^{\lambda }
\end{equation*}%
for every $\tau \in \left[ 0,T\right] .$

(2) Whereas the Bernstein process of Proposition 1 is stationary, that of
Proposition 2 is not. This is intuitively understandable, as some kind of
non trivial dynamics ought to be necessary to steer the process from a
deterministic state at $t=0$ to a Gaussian distribution at the end of its
journey.

\bigskip

Of course, there is an almost unlimited number of possibilities of getting
various Bernstein processes in the above manner, just by choosing $\varphi
_{0,\lambda }$ and $\psi _{T,\lambda }$ appropriately. We complete this
section by constructing yet another process which shares many properties of
a Markovian bridge. For this, we consider the forward-backward system (\ref%
{forwardharmonic})-(\ref{backwardharmonic}) with Dirac measures concentrated
at the origin and at a given point $\mathsf{a}\in \mathbb{R}^{N}$ as
initial-final data, namely, 
\begin{align}
\partial _{t}u(\mathsf{x},t)& =\frac{1}{2}\Delta _{\mathsf{x}}u(\mathsf{x}%
,t)-\frac{\lambda ^{2}\left\vert \mathsf{x}\right\vert ^{2}}{2}u(\mathsf{x}%
,t),\text{ \ \ }(\mathsf{x},t)\in \mathbb{R}^{N}\mathbb{\times }\left( 0,T%
\right] ,  \notag \\
u(\mathsf{x},0)& =\varphi _{0,\lambda }(\mathsf{x})=\mathsf{m}_{\lambda
}\delta (\mathsf{x})  \label{initialvalueproblem}
\end{align}%
and%
\begin{align}
-\partial _{t}v(\mathsf{x},t)& =\frac{1}{2}\Delta _{\mathsf{x}}v(\mathsf{x}%
,t)-\frac{\lambda ^{2}\left\vert \mathsf{x}\right\vert ^{2}}{2}v(\mathsf{x}%
,t),\text{ \ \ }(\mathsf{x},t)\in \mathbb{R}^{N}\mathbb{\times }\left[
0,T\right) ,  \notag \\
v(\mathsf{x},T)& =\psi _{T,\lambda }(\mathsf{x})=\mathsf{m}_{\lambda }\delta
(\mathsf{x}-\mathsf{a}),  \label{finalvalueproblem}
\end{align}%
where 
\begin{equation*}
\mathsf{m}_{\lambda }=\left( \frac{2\pi \sinh \left( \lambda T\right) }{%
\lambda }\right) ^{\frac{N}{4}}\exp \left[ \frac{\lambda \coth \left(
\lambda T\right) \left\vert \mathsf{a}\right\vert ^{2}}{4}\right] .
\end{equation*}%
It is easily verified that the normalization condition%
\begin{equation*}
\mathsf{m}_{\lambda }^{2}\int_{\mathbb{R}^{N}\times \mathbb{R}^{N}}d\mathsf{x%
}d\mathsf{y}\delta (\mathsf{x})\delta (\mathsf{y-a})g_{\lambda }(\mathsf{x}%
,T,\mathsf{y})=1
\end{equation*}%
holds with Mehler's kernel given above, and that the unique classical
positive solutions to (\ref{initialvalueproblem})-(\ref{finalvalueproblem})
which satisfy (\ref{requirement}) are now%
\begin{equation}
u_{\lambda }(\mathsf{x},t)=\mathsf{n}_{\lambda }\sinh ^{-\frac{N}{2}}\left(
\lambda t\right) \exp \left[ -\frac{\alpha _{\lambda }(t)\left\vert \mathsf{x%
}\right\vert ^{2}}{2}\right]  \label{forwardsolutionter}
\end{equation}%
and%
\begin{eqnarray}
v_{\lambda }(\mathsf{x},t) &=&\mathsf{n}_{\lambda }\exp \left[ -\frac{\alpha
_{\lambda }(T-t)\left\vert \mathsf{a}\right\vert ^{2}}{2}\right] \sinh ^{-%
\frac{N}{2}}\left( \lambda \left( T-t\right) \right)
\label{backwardsolutionter} \\
&&\times \exp \left[ -\frac{1}{2}\left( \alpha _{\lambda }(T-t)\left\vert 
\mathsf{x}\right\vert ^{2}-\frac{2\lambda \left( \mathsf{a},\mathsf{x}%
\right) _{\mathbb{R}^{N}}}{\sinh \left( \lambda \left( T-t\right) \right) }%
\right) \right] ,
\end{eqnarray}%
respectively, where%
\begin{equation}
\alpha _{\lambda }(t)=\lambda \coth \left( \lambda t\right)  \label{alpha}
\end{equation}%
and%
\begin{equation}
\mathsf{n}_{\lambda }=\left( \frac{\lambda \sinh \left( \lambda T\right) }{%
2\pi }\right) ^{\frac{N}{4}}\exp \left[ \frac{\alpha _{\lambda
}(T)\left\vert \mathsf{a}\right\vert ^{2}}{4}\right] .  \label{coefficient}
\end{equation}%
Then the following result is valid:

\bigskip

\textbf{Proposition 3.}\textit{\ The Bernstein process }$Z_{\tau \in \left[
0,T\right] }^{\lambda }$\textit{\ associated with (\ref{initialvalueproblem}%
)-(\ref{finalvalueproblem}) in the sense of Corollary 1 is a Gaussian and
Markovian process such that}%
\begin{equation}
\mathbb{P}_{\mu }\left( Z_{t}^{\lambda }\in E\right) =(2\pi \sigma _{\lambda
}(t))^{-\frac{N}{2}}\int_{E}d\mathsf{x}\exp \left[ -\frac{\left\vert \mathsf{%
x-a}_{\lambda }(t)\right\vert ^{2}}{2\sigma _{\lambda }(t)}\right]
\label{gaussianprocessbis}
\end{equation}%
\textit{for each }$t\in \left( 0,T\right) $\textit{\ and every} $E\in 
\mathcal{B}_{N}$, \textit{where}%
\begin{equation}
\mathsf{a}_{\lambda }(t)=\frac{\sinh (\lambda t)}{\sinh (\lambda T)}\mathsf{a%
\label{meanvector}}
\end{equation}%
\textit{and}%
\begin{equation}
\sigma _{\lambda }(t)=\frac{\sinh \left( \lambda (T-t\right) )\sinh (\lambda
t)}{\lambda \sinh (\lambda T)}.  \label{variance}
\end{equation}%
\textit{Furthermore we have}%
\begin{eqnarray}
&&\mathbb{P}_{\mu }\left( Z_{0}^{\lambda }=\mathsf{o}\right)  \notag \\
&=&\mathbb{P}_{\mu }\left( Z_{T}^{\lambda }=\mathsf{a}\right) =1,
\label{probability7}
\end{eqnarray}%
\textit{and the components of }$Z_{\tau \in \left[ 0,T\right] }^{\lambda }$%
\textit{\ satisfy the relation}%
\begin{equation}
\mathbb{E}_{\mu }\left( (Z_{s}^{\lambda ,i}-a_{\lambda
}^{i}(s))(Z_{t}^{\lambda ,j}-a_{\lambda }^{j}(t))\right) =\left\{ 
\begin{array}{c}
\frac{\sinh \left( \lambda (T-t\right) )\sinh (\lambda s)}{\lambda \sinh
(\lambda T)}\delta _{i,j}\text{ \ \ \ for }t\geq s, \\ 
\\ 
\frac{\sinh \left( \lambda (T-s\right) )\sinh (\lambda t)}{\lambda \sinh
(\lambda T)}\delta _{i,j}\text{ \ \ for }t\leq s%
\end{array}%
\right.  \label{covariance}
\end{equation}%
\textit{for all }$s,t\in \left[ 0,T\right] $\textit{\ and all }$i,j\in
\left\{ 1,...,N\right\} $\textit{. In fact, }$Z_{\tau \in \left[ 0,T\right]
}^{\lambda }$\textit{\ is a non-stationary process pinned down at the origin
when }$t=0$, \textit{at} $\mathsf{a}$ \textit{when} $t=T$, \textit{and
exhibiting maximal randomness when }$t=\frac{T}{2}.$

\bigskip

\textbf{Proof.} We begin by proving (\ref{gaussianprocessbis}). Using (\ref%
{forwardsolutionter}), (\ref{backwardsolutionter}) and (\ref{coefficient})
we first have%
\begin{eqnarray}
&&u_{\lambda }(\mathsf{x},t)v_{\lambda }(\mathsf{x},t)  \notag \\
&=&\left( \frac{\lambda \sinh (\lambda T)}{2\pi \sinh \left( \lambda
(T-t\right) )\sinh (\lambda t)}\right) ^{\frac{N}{2}}\exp \left[ \frac{%
\left( \alpha _{\lambda }(T)-\alpha _{\lambda }(T-t)\right) \left\vert 
\mathsf{a}\right\vert ^{2}}{2}\right]  \notag \\
&&\times \exp \left[ -\frac{1}{2}\left( \left( \alpha _{\lambda }(t)+\alpha
_{\lambda }(T-t)\right) \left\vert \mathsf{x}\right\vert ^{2}-\frac{2\lambda
\left( \mathsf{a},\mathsf{x}\right) _{\mathbb{R}^{N}}}{\sinh \left( \lambda
\left( T-t\right) \right) }\right) \right]  \label{probadensity}
\end{eqnarray}%
after regrouping terms, and furthermore%
\begin{eqnarray}
\alpha _{\lambda }(T)-\alpha _{\lambda }(T-t) &=&-\frac{\lambda \sinh
(\lambda t)}{\sinh \left( \lambda (T-t\right) )\sinh (\lambda T)}
\label{coefficient2} \\
\alpha _{\lambda }(t)+\alpha _{\lambda }(T-t) &=&\frac{\lambda \sinh
(\lambda T)}{\sinh \left( \lambda (T-t\right) )\sinh (\lambda t)}
\label{coefficient3}
\end{eqnarray}%
from (\ref{alpha}). The substitution of (\ref{coefficient2})-(\ref%
{coefficient3}) into (\ref{probadensity}) then leads to%
\begin{eqnarray}
&&u_{\lambda }(\mathsf{x},t)v_{\lambda }(\mathsf{x},t)  \notag \\
&=&\left( \frac{\lambda \sinh (\lambda T)}{2\pi \sinh \left( \lambda
(T-t\right) )\sinh (\lambda t)}\right) ^{\frac{N}{2}}\exp \left[ -\frac{%
\lambda \sinh (\lambda t)\left\vert \mathsf{a}\right\vert ^{2}}{2\sinh
\left( \lambda (T-t\right) )\sinh (\lambda T)}\right]  \notag \\
&&\times \exp \left[ -\frac{\lambda }{2}\left( \frac{\sinh (\lambda
T)\left\vert \mathsf{x}\right\vert ^{2}-2\sinh (\lambda t)\left( \mathsf{a},%
\mathsf{x}\right) _{\mathbb{R}^{N}}}{\sinh \left( \lambda (T-t\right) )\sinh
(\lambda t)}\right) \right] .  \label{probadensitybis}
\end{eqnarray}%
Now, for the numerator of the argument in the second exponential of (\ref%
{probadensitybis}) we have 
\begin{eqnarray}
&&\sinh (\lambda T)\left\vert \mathsf{x}\right\vert ^{2}-2\sinh (\lambda
t)\left( \mathsf{a},\mathsf{x}\right) _{\mathbb{R}^{N}}  \notag \\
&=&\sinh (\lambda T)\left\vert \mathsf{x-a}_{\lambda }(t)\right\vert ^{2}-%
\frac{\sinh ^{2}(\lambda t)\left\vert \mathsf{a}\right\vert ^{2}}{\sinh
(\lambda T)}  \label{identity}
\end{eqnarray}%
by virtue of (\ref{meanvector}). Therefore, taking (\ref{variance}) and (\ref%
{identity}) into account in (\ref{probadensitybis}) we get%
\begin{eqnarray*}
&&u_{\lambda }(\mathsf{x},t)v_{\lambda }(\mathsf{x},t) \\
&=&(2\pi \sigma _{\lambda }(t))^{-\frac{N}{2}}\exp \left[ -\frac{\left\vert 
\mathsf{x-a}_{\lambda }(t)\right\vert ^{2}}{2\sigma _{\lambda }(t)}\right]
\end{eqnarray*}%
following the cancellation of two exponential factors, which gives the
desired result according to (\ref{probability4}).

As for the proof of (\ref{probability7}), we remark that (\ref{probability5}%
) and (\ref{probability6}) lead to%
\begin{eqnarray*}
&&\mathbb{P}_{\mu }\left( Z_{0}\in E\right) \\
&=&\int_{E}d\mathsf{x}\delta (\mathsf{x})\exp \left[ -\frac{\lambda \left(
\cosh (\lambda T)\left\vert \mathsf{x}\right\vert ^{2}-2\left( \mathsf{a,x}%
\right) _{\mathbb{R}^{N}}\right) }{2\sinh \left( \lambda T\right) }\right]
\end{eqnarray*}%
and%
\begin{eqnarray*}
&&\mathbb{P}_{\mu }\left( Z_{T}\in E\right) \\
&=&\exp \left[ \frac{\alpha _{\lambda }(T)\left\vert \mathsf{a}\right\vert
^{2}}{2}\right] \int_{E}d\mathsf{x}\delta (\mathsf{x-a})\exp \left[ -\frac{%
\alpha _{\lambda }(T)\left\vert \mathsf{x}\right\vert ^{2}}{2}\right]
\end{eqnarray*}%
for every $E\in \mathcal{B}_{N}$, respectively, which immediately imply the
claim.

We now turn to the proof of (\ref{covariance}), by noticing that in this
case the finite-dimensional density in $\mathbb{R}^{nN}$ is%
\begin{eqnarray*}
&&\dprod\limits_{k=2}^{n}g_{\lambda }\left( \mathsf{x}_{k},t_{k}-t_{k-1},%
\mathsf{x}_{k-1}\right) \times u_{\lambda }(\mathsf{x}_{1},t_{1})v_{\lambda
}(\mathsf{x}_{n},t_{n}) \\
&=&\left( 2\pi \right) ^{-\frac{nN}{2}}\left( \frac{\lambda ^{n}\sinh
(\lambda T)}{\sinh \left( \lambda (T-t_{n}\right) )\sinh (\lambda t_{1})}%
\right) ^{\frac{N}{2}}\left( \dprod\limits_{k=2}^{n}\sinh \left( \lambda
(t_{k}-t_{k-1})\right) \right) ^{-\frac{N}{2}} \\
&&\times \exp \left[ \frac{1}{2}(\alpha _{\lambda }(T)-\alpha _{\lambda
}(T-t_{n}))\left\vert \mathsf{a}\right\vert ^{2}\right] \\
&&\times \exp \left[ -\frac{\lambda }{2}\sum_{k=2}^{n}\frac{\cosh (\lambda
(t_{k}-t_{k-1}))\left( \left\vert \mathsf{x}_{k}\right\vert ^{2}+\left\vert 
\mathsf{x}_{k-1}\right\vert ^{2}\right) -2\left( \mathsf{x}_{k}\mathsf{,x}%
_{k-1}\right) _{\mathbb{R}^{N}}}{\sinh \left( \lambda (t_{k}-t_{k-1})\right) 
}\right] \\
&&\times \exp \left[ -\frac{1}{2}\left( \alpha _{\lambda }(t_{1})\left\vert 
\mathsf{x}_{1}\right\vert ^{2}+\alpha _{\lambda }(T-t_{n})\left\vert \mathsf{%
x}_{n}\right\vert ^{2}\right) \right] \times \exp \left[ \frac{\lambda
\left( \mathsf{a},\mathsf{x}_{n}\right) _{\mathbb{R}^{N}}}{\sinh \left(
\lambda (T-t_{n}\right) )}\right] .
\end{eqnarray*}%
Therefore, for the tridiagonal matrix $C_{\lambda }^{-1}$ corresponding to
the quadratic part when $N=1$ we obtain%
\begin{equation*}
C_{\lambda ,k,k}^{-1}=\left\{ 
\begin{array}{c}
\frac{\lambda \sinh (\lambda t_{2})}{\sinh \left( \lambda
(t_{2}-t_{1}\right) )\sinh (\lambda t_{1})}\text{ \ \ for }k=1, \\ 
\\ 
\frac{\lambda \sinh (\lambda (t_{k+1}-t_{k-1}))}{\sinh \left( \lambda
(t_{k+1}-t_{k}\right) )\sinh (\lambda (t_{k}-t_{k-1}))}\text{ \ \ for }%
k=2,...,n-1, \\ 
\\ 
\frac{\lambda \sinh (\lambda (T-t_{n-1}))}{\sinh \left( \lambda
(T-t_{n}\right) )\sinh (\lambda (t_{n}-t_{n-1}))}\text{ \ \ for }k=n%
\end{array}%
\right.
\end{equation*}%
(the second line still not being there if $n=2)$, and 
\begin{equation*}
C_{\lambda ,k,k-1}^{-1}=C_{\lambda ,k-1,k}^{-1}=-\frac{\lambda }{\sinh
\left( \lambda \left\vert t_{k}-t_{k-1}\right\vert \right) }\text{ \ \ for }%
k=2,...,n.
\end{equation*}%
Consequently, inverting again and using numerous relations among hyperbolic
functions we eventually get%
\begin{equation*}
C_{\lambda ,k,l}=\left\{ 
\begin{array}{c}
\frac{\sinh \left( \lambda (T-t_{k})\right) \sinh \left( \lambda
t_{l}\right) }{\lambda \sinh \left( \lambda T\right) }\text{ \ \ for }k\geq
l, \\ 
\\ 
\frac{\sinh \left( \lambda (T-t_{l})\right) \sinh \left( \lambda
t_{k}\right) }{\lambda \sinh \left( \lambda T\right) }\text{ \ \ for }k\leq
l,%
\end{array}%
\right.
\end{equation*}%
which leads to (\ref{covariance}) by standard arguments. Finally, we note
that the curve $\sigma _{\lambda }:\left[ 0,T\right] \mapsto \mathbb{R}%
_{0}^{+}$ given by (\ref{variance}) is concave aside from satisfying $\sigma
_{\lambda }(0)=\sigma _{\lambda }(T)=0$, and that it takes on the maximal
value%
\begin{equation*}
\sigma _{\lambda }\left( \frac{T}{2}\right) =\frac{\sinh ^{2}\left( \frac{%
\lambda T}{2}\right) }{\lambda \sinh (\lambda T)}.
\end{equation*}%
Thus, the process $Z_{\tau \in \left[ 0,T\right] }^{\lambda }$ is indeed
non-stationary with the stated properties. \ \ $\blacksquare $

\bigskip

\textsc{Remark.} We may dub the process $Z_{\tau \in \left[ 0,T\right]
}^{\lambda }$ of the preceding corollary a \textit{Bernstein bridge}, which
represents a random curve whose ends are pinned down at specified points in
space. We remark that the corresponding Gaussian law is no longer centered,
unless $\mathsf{a}=\mathsf{o}$ in which case the process materializes a
Markovian loop which retains the main features of a Brownian loop. In fact, $%
Z_{\tau \in \left[ 0,T\right] }^{\lambda }$ does reduce to a Brownian bridge
in the limit $\lambda \rightarrow 0_{+}$ since%
\begin{equation*}
\lim_{\lambda \rightarrow 0_{+}}\mathbb{E}_{\mu }\left( (Z_{s}^{\lambda
,i}-a_{\lambda }^{i}(s))(Z_{t}^{\lambda ,j}-a_{\lambda }^{j}(t))\right)
=\left\{ 
\begin{array}{c}
\frac{(T-t)s}{T}\delta _{i,j}\text{ \ \ \ for }t\geq s, \\ 
\\ 
\frac{(T-s)t}{T}\delta _{i,j}\text{ \ \ for }t\leq s%
\end{array}%
\right.
\end{equation*}%
according to (\ref{covariance}).

\bigskip

In the next section we introduce a new class of Bernstein processes which we
can eventually relate to the so-called periodic Ornstein-Uhlenbeck process,
and which is generated by a one-parameter family of non-Markovian
probability measures.

\section{A family of non-Markovian Bernstein processes and the periodic
Ornstein-Uhlenbeck process}

The Bernstein processes of this section are still defined from measures
which are intimately tied up with problems of the form (\ref{forwardharmonic}%
)-(\ref{backwardharmonic}), but their finite-dimensional distributions will
be determined exclusively from Theorem 1. Instead of considering just one
pair of parabolic problems such as (\ref{forwardharmonic})-(\ref%
{backwardharmonic}), we first introduce an infinite hierarchy of
forward-backward systems of the form%
\begin{align}
\partial _{t}u(\mathsf{x},t)& =\frac{1}{2}\Delta _{\mathsf{x}}u(\mathsf{x}%
,t)-\frac{\lambda ^{2}}{2}\left\vert \mathsf{x}\right\vert ^{2}u(\mathsf{x}%
,t),\text{ \ }(\mathsf{x},t)\in \mathbb{R}^{N}\times \left( 0,T\right] , 
\notag \\
u(\mathsf{x},0)& =\varphi _{\mathsf{m},0,\lambda }(\mathsf{x})=e^{\frac{1}{2}%
\left( \sum_{j=1}^{N}m_{j}+\frac{N}{2}\right) \lambda T}\otimes
_{j=1}^{N}h_{m_{j},\lambda }(\mathsf{x}),\text{\ \ }\mathsf{x}\in \mathbb{R}%
^{N}  \label{forwardharmonicter}
\end{align}%
and%
\begin{align}
-\partial _{t}v(\mathsf{x},t)& =\frac{1}{2}\Delta _{\mathsf{x}}v(\mathsf{x}%
,t)-\frac{\lambda ^{2}}{2}\left\vert \mathsf{x}\right\vert ^{2}v(\mathsf{x}%
,t),\text{ \ }(\mathsf{x},t)\in \mathbb{R}^{N}\times \left[ 0,T\right) , 
\notag \\
v(\mathsf{x,}T)& =\psi _{\mathsf{m},T,\lambda }(\mathsf{x})=e^{\frac{1}{2}%
\left( \sum_{j=1}^{N}m_{j}+\frac{N}{2}\right) \lambda T}\otimes
_{j=1}^{N}h_{m_{j},\lambda }(\mathsf{x}),\text{ \ }\mathsf{x}\in \mathbb{R}%
^{N},  \label{backwardharmonicter}
\end{align}%
namely, one such pair for each $m_{j}\in \mathbb{N}$ and every $j$, where
the $h_{m_{j},\lambda }$'s are the Hermite functions of Lemma 1.
Accordingly, this means that we are considering as many pairs of such
systems as is necessary to take into account the whole pure point spectrum
of the elliptic operator on the right-hand side. We remark that (\ref%
{forwardharmonicter})-(\ref{backwardharmonicter}) constitutes a
generalization of (\ref{initial2})-(\ref{final2}), the latter system being
associated with the bottom of the spectrum where $m_{j}\mathsf{=0}$ for each 
$j$. Whereas the associated measures remain suitably normalized, the
drawback is that they are no longer positive according to the following
result:

\bigskip

\textbf{Lemma 2.} \textit{Let us write }$\mathsf{m:=(}m_{1},...,m_{N})\in 
\mathbb{N}^{N}$ \textit{and let us consider the sequence of measures}%
\begin{equation}
\mu _{\mathsf{m},\lambda }\left( F\right) =\int_{F}d\mu _{\mathsf{m},\lambda
}\left( \mathsf{x},\mathsf{y}\right)  \label{signedmeasures}
\end{equation}%
\textit{where}%
\begin{equation}
\mu _{\mathsf{m},\lambda }\left( \mathsf{x},\mathsf{y}\right) =\varphi _{%
\mathsf{m},0,\lambda }\left( \mathsf{x}\right) \psi _{\mathsf{m},T,\lambda
}\left( \mathsf{y}\right) g_{\lambda }(\mathsf{x},T,\mathsf{y})
\label{densitybis}
\end{equation}%
\textit{\ and }$F\in \mathcal{B}_{N}\times \mathcal{B}_{N}$\textit{. Then we
have}%
\begin{equation}
\int_{\mathbb{R}^{N}\times \mathbb{R}^{N}}d\mu _{\mathsf{m},\lambda }\left( 
\mathsf{x},\mathsf{y}\right) =1  \label{normalizationter}
\end{equation}%
\textit{and for every} $\mathsf{m\neq 0}$ \textit{the} $\mu _{\mathsf{m}%
,\lambda }$\textit{'s are signed measures.}

\bigskip

\textbf{Proof. }According to Proposition A.1 of the appendix we have%
\begin{eqnarray}
&&g_{\lambda }(\mathsf{x},T,\mathsf{y})  \notag \\
&=&\sum_{\mathsf{n\in }\mathbb{N}^{N}}e^{-\left( \sum_{j=1}^{N}n_{j}+\frac{N%
}{2}\right) \lambda T}\otimes _{j=1}^{N}h_{n_{j},\lambda }\left( \mathsf{x}%
\right) \times \otimes _{j=1}^{N}h_{n_{j},\lambda }\left( \mathsf{y}\right)
\label{expansionbis}
\end{eqnarray}%
for Mehler's kernel (\ref{Mehlerkernel}), where the series converges
absolutely and uniformly for all $\mathsf{x},\mathsf{y}$\textsf{\ }$\in 
\mathbb{R}^{N}$. Moreover, since the tensor products $\otimes
_{j=1}^{N}h_{m_{j},\lambda }$ provide an\textit{\ }orthonormal basis of $%
L^{2}\left( \mathbb{R}^{N},\mathbb{C}\right) $ we have%
\begin{equation}
\int_{\mathbb{R}^{N}}d\mathsf{x}\otimes _{j=1}^{N}h_{m_{j},\lambda }\left( 
\mathsf{x}\right) \times \otimes _{j=1}^{N}h_{n_{j},\lambda }\left( \mathsf{x%
}\right) =\dprod\limits_{j=1}^{N}\delta _{m_{j},n_{j}}  \label{orthogonality}
\end{equation}%
for all $m_{j},n_{j}\in \mathbb{N}$. Consequently, substituting $\varphi _{%
\mathsf{m},0,\lambda }$\textsf{, }$\psi _{\mathsf{m},T,\lambda }$\textsf{\ }%
and\textsf{\ }(\ref{expansionbis}) into the left-hand side of (\ref%
{normalizationter}) and taking (\ref{orthogonality}) into account we obtain%
\begin{eqnarray*}
&&\int_{\mathbb{R}^{N}\times \mathbb{R}^{N}}d\mathsf{x}d\mathsf{y}\varphi _{%
\mathsf{m},0,\lambda }\left( \mathsf{x}\right) \psi _{\mathsf{m},T,\lambda
}\left( \mathsf{y}\right) g_{\lambda }(\mathsf{x},T,\mathsf{y}) \\
&=&\sum_{\mathsf{n\in }\mathbb{N}^{N}}e^{\sum_{j=1}^{N}(m_{j}-n_{j})\lambda
T}\dprod\limits_{j=1}^{N}\delta _{m_{j},n_{j}}=1.
\end{eqnarray*}%
Finally, the $\mu _{\mathsf{m},\lambda }$'s are signed measures because of
the existence of real zeroes for each $h_{m_{j},\lambda }$ when $m_{j}\neq 0$%
. \ \ $\blacksquare $

\bigskip

The fact that the above measures are signed prevents one from applying
directly the general results of Section 2 to (\ref{forwardharmonicter})-(\ref%
{backwardharmonicter}) with a fixed $\mathsf{m\neq 0}$, as it would prevent
one from applying the main result of \cite{beurling} briefly discussed at
the end of Section 2. However, we can still save the day by constructing a
one-parameter family of \textit{bona fide} probability measures from all the 
$\mu _{\mathsf{m},\lambda }$'s. Indeed the following result is valid, where
we call a measure non-Markovian whenever (\ref{markovianmeasures}) does 
\textit{not} hold:

\bigskip

\textbf{Lemma 3.} \textit{For each }$\lambda >0$ \textit{there exists a
one-parameter family }$\left( \hat{\mu}_{\lambda ,\theta }\right) _{\theta
>0}$\textit{\ of positive, non-Markovian measures of the form}%
\begin{equation*}
\hat{\mu}_{\lambda ,\theta }(F)=\int_{F}d\hat{\mu}_{\lambda ,\theta }\left( 
\mathsf{x,y}\right)
\end{equation*}%
\textit{where} $F\in \mathcal{B}_{N}\times \mathcal{B}_{N}$\textit{, which
satisfy} 
\begin{equation}
\int_{\mathbb{R}^{N}\times \mathbb{R}^{N}}d\hat{\mu}_{\lambda ,\theta
}\left( \mathsf{x,y}\right) =1  \label{normalizationquarto}
\end{equation}%
\textit{and which disintegrate into a statistical mixture of the }$\mu _{%
\mathsf{m},\lambda }$\textit{'s. In other words, for each }$\lambda >0$, 
\textit{each }$\theta >0$ \textit{and every }$\mathsf{m}\in \mathbb{N}^{N}$%
\textit{\ there exist numbers }$p_{\mathsf{m},\lambda ,\theta }>0$\textit{\
such that} 
\begin{equation*}
\hat{\mu}_{\lambda ,\theta }=\sum_{\mathsf{m\in }\mathbb{N}^{N}}p_{\mathsf{m}%
,\lambda ,\theta }\mu _{\mathsf{m,}\lambda }\text{ \ \textit{where }}\sum_{%
\mathsf{m\in }\mathbb{N}^{N}}p_{\mathsf{m},\lambda ,\theta }=1\text{.}
\end{equation*}%
\textit{In fact, it is sufficient to take}%
\begin{equation}
\hat{\mu}_{\lambda ,\theta }\left( \mathsf{x,y}\right) =\left( 2\left( \cosh
(\lambda \left( \theta +1\right) T)-1\right) \right) ^{\frac{N}{2}%
}g_{\lambda }(\mathsf{x},T,\mathsf{y})g_{\lambda }(\mathsf{x},\theta T,%
\mathsf{y}).  \label{density}
\end{equation}

\bigskip

\textbf{Proof.} From the series expansion (\ref{expansionbis}) and the
orthogonality relations (\ref{orthogonality}) we have%
\begin{eqnarray*}
&&\int_{\mathbb{R}^{N}\times \mathbb{R}^{N}}d\mathsf{x}d\mathsf{y}g_{\lambda
}(\mathsf{x},T,\mathsf{y})g_{\lambda }(\mathsf{x},\theta T,\mathsf{y}) \\
&=&\sum_{\mathsf{m\in }\mathbb{N}^{N}}e^{-\left( \sum_{j=1}^{N}m_{j}+\frac{N%
}{2}\right) \lambda \left( \theta +1\right) T} \\
&=&\left( 2\sinh \left( \lambda \left( \theta +1\right) \frac{T}{2}\right)
\right) ^{-N}
\end{eqnarray*}%
by summing the underlying geometric series, from which we obtain (\ref%
{normalizationquarto}) as a consequence of (\ref{density}) and the identity%
\begin{equation*}
4\sinh ^{2}\left( \lambda \left( \theta +1\right) \frac{T}{2}\right)
=2\left( \cosh (\lambda \left( \theta +1\right) T)-1\right) .
\end{equation*}%
As for the second statement of the lemma, let us define the sequence of
numbers%
\begin{equation}
p_{\mathsf{m},\lambda ,\theta }:=\left( 2\left( \cosh (\lambda \left( \theta
+1\right) T)-1\right) \right) ^{\frac{N}{2}}e^{-\left( \sum_{j=1}^{N}m_{j}+%
\frac{N}{2}\right) \lambda \left( \theta +1\right) T}>0.  \label{sequence}
\end{equation}%
Summing as above we get%
\begin{equation*}
\sum_{\mathsf{m\in }\mathbb{N}^{N}}p_{\mathsf{m},\lambda ,\theta }=1
\end{equation*}%
as required. Therefore, taking into account (\ref{densitybis}), (\ref%
{expansionbis}) and (\ref{sequence}) we obtain%
\begin{eqnarray*}
&&\sum_{\mathsf{m\in }\mathbb{N}^{N}}p_{\mathsf{m},\lambda ,\theta }\mu _{%
\mathsf{m},\lambda }\left( \mathsf{x},\mathsf{y}\right) \\
&=&g_{\lambda }(\mathsf{x},T,\mathsf{y})\sum_{\mathsf{m\in }\mathbb{N}%
^{N}}p_{\mathsf{m},\lambda ,\theta }\varphi _{\mathsf{m},0,\lambda }\left( 
\mathsf{x}\right) \psi _{\mathsf{m},T,\lambda }\left( \mathsf{y}\right) \\
&=&\left( 2\left( \cosh (\lambda \left( \theta +1\right) T)-1\right) \right)
^{\frac{N}{2}}g_{\lambda }(\mathsf{x},T,\mathsf{y}) \\
&&\times \sum_{\mathsf{m\in }\mathbb{N}^{N}}e^{-\left( \sum_{j=1}^{N}m_{j}+%
\frac{N}{2}\right) \lambda \theta T}\otimes _{j=1}^{N}h_{m_{j},\lambda
}\left( \mathsf{x}\right) \times \otimes _{j=1}^{N}h_{m_{j},\lambda }\left( 
\mathsf{y}\right) \\
&=&\hat{\mu}_{\lambda ,\theta }\left( \mathsf{x,y}\right)
\end{eqnarray*}%
according to (\ref{density}), which is the desired conclusion. \ \ $%
\blacksquare $

\bigskip

\textsc{Remark.} Strictly speaking, the measure $\hat{\mu}_{\lambda ,\theta
} $ does not exist for $\theta =0$ but the limit%
\begin{equation}
\hat{\mu}_{\lambda ,+}\left( \mathsf{x,y}\right) :=\lim_{\theta \rightarrow
0_{+}}\hat{\mu}_{\lambda ,\theta }\left( \mathsf{x,y}\right) =\left( 2\left(
\cosh (\lambda T)-1\right) \right) ^{\frac{N}{2}}g_{\lambda }(\mathsf{x},T,%
\mathsf{y})\delta (\mathsf{x}-\mathsf{y})  \label{limitingcase}
\end{equation}%
does, by virtue of the fact that $g_{\lambda }$ is Green's function
associated with the partial differential equation in (\ref{initial2}). Said
differently, Lemma 3 and its proof remain valid for the measure $\hat{\mu}%
_{\lambda ,+}$ associated with (\ref{limitingcase}) since we have%
\begin{equation*}
\sum_{\mathsf{m\in }\mathbb{N}^{N}}\otimes _{j=1}^{N}h_{m_{j},\lambda
}\left( \mathsf{x}\right) \times \otimes _{j=1}^{N}h_{m_{j},\lambda }\left( 
\mathsf{y}\right) =\delta (\mathsf{x}-\mathsf{y})
\end{equation*}%
in the sense of distributions for all $\mathsf{x},\mathsf{y\in }\mathbb{R}%
^{N}$, as a consequence of the completeness of the $\otimes
_{j=1}^{N}h_{m_{j},\lambda }$'s in $L^{2}\left( \mathbb{R}^{N},\mathbb{C}%
\right) $.

\bigskip

The above developments now lead to the following result:

\bigskip

\textbf{Proposition 4. }\textit{For each }$\lambda >0$, \textit{the
Bernstein processes }$Z_{\tau \in \left[ 0,T\right] }^{\lambda ,\theta >0}$%
\textit{\ associated with the above infinite hierarchy through the measures }%
$\hat{\mu}_{\lambda ,\theta }$ \textit{are stationary Gaussian and
non-Markovian processes such that}%
\begin{equation}
\mathbb{P}_{\mu }\left( Z_{t}^{\lambda ,\theta }\in E\right) =(2\pi \sigma
_{\lambda ,\theta })^{-\frac{N}{2}}\int_{E}dx\exp \left[ -\frac{\left\vert 
\mathsf{x}\right\vert ^{2}}{2\sigma _{\lambda ,\theta }}\right]
\label{gaussianprocesster}
\end{equation}%
\textit{for each }$t\in \left[ 0,T\right] $ \textit{and every} $E\in 
\mathcal{B}_{N}$, \textit{where}%
\begin{equation}
\sigma _{\lambda ,\theta }=\frac{\sinh \left( \lambda \left( \theta
+1\right) T\right) }{2\lambda \left( \cosh (\lambda \left( \theta +1\right)
T)-1\right) }.  \label{variancebis}
\end{equation}%
\textit{Furthermore, the components of }$Z_{\tau \in \left[ 0,T\right]
}^{\lambda ,\theta }$\textit{\ satisfy the relation}%
\begin{equation}
\mathbb{E}_{\mu }\left( Z_{s}^{\lambda ,\theta ,i}Z_{t}^{\lambda ,\theta
,j}\right) =\frac{\cosh \left( \lambda \left( \left\vert t-s\right\vert -%
\frac{\left( \theta +1\right) T}{2}\right) \right) }{2\lambda \sinh \left( 
\frac{\lambda \left( \theta +1\right) T}{2}\right) }\delta _{i,j}
\label{covariancebis}
\end{equation}%
\textit{for all }$s,t\in \left[ 0,T\right] $\textit{\ and all }$i,j\in
\left\{ 1,...,N\right\} $\textit{. Finally, the process }$Z_{\tau \in \left[
0,T\right] }^{\lambda ,+}$\textit{\ associated with the measure }$\hat{\mu}%
_{\lambda ,+}$ \textit{identifies in law with a stationary process} \textit{%
whose components are all independent, one-dimensional and periodic
Ornstein-Uhlenbeck processes.}

\bigskip

\textbf{Proof.} The processes $Z_{\tau \in \left[ 0,T\right] }^{\lambda
,\theta >0}$ are Gaussian and non-Markovian by virtue of (\ref{distribution}%
) with Green's function (\ref{Mehlerkernel}) and the measures $\hat{\mu}%
_{\lambda ,\theta }$. Furthermore, using the symmetry properties of $%
g_{\lambda }$ with respect to the spatial variables, and twice the related
semi-group composition law (\ref{semigroup}) of the appendix, (\ref%
{probability1}) with (\ref{density}) becomes 
\begin{eqnarray*}
&&\mathbb{P}_{\mu }\left( Z_{t}^{\lambda ,\theta }\in E\right) \\
&=&\left( 2\left( \cosh (\lambda \left( \theta +1\right) T)-1\right) \right)
^{\frac{N}{2}}\int_{E}d\mathsf{x}g_{\lambda }\left( \mathsf{x},\left( \theta
+1\right) T,\mathsf{x}\right)
\end{eqnarray*}%
for every $t\in (0,T)$. The same result obtains for $t=0$ and $t=T$ as\ a
consequence of (\ref{probability2})-(\ref{probability3}), so that (\ref%
{gaussianprocesster})-(\ref{variancebis}) follows immediately since%
\begin{eqnarray*}
&&g_{\lambda }(\mathsf{x},\left( \theta +1\right) T,\mathsf{x}) \\
&=&\left( \frac{\lambda }{2\pi \sinh \left( \lambda \left( \theta +1\right)
T\right) }\right) ^{\frac{N}{2}}\exp \left[ -\frac{\lambda \left( \cosh
(\lambda \left( \theta +1\right) T)-\mathsf{1}\right) \left\vert \mathsf{x}%
\right\vert ^{2}}{\sinh \left( \lambda \left( \theta +1\right) T\right) }%
\right] .
\end{eqnarray*}%
We now turn to the proof of (\ref{covariancebis}) by determining the
Gaussian density of $\left( Z_{t_{1}}^{\lambda ,\theta
},...,Z_{t_{n}}^{\lambda ,\theta }\right) $. Thus, using $\hat{\mu}_{\lambda
,\theta }$ in (\ref{distribution}) and integrating first over $\mathsf{y}$
and then over $\mathsf{x}$ we may write%
\begin{eqnarray*}
&&\mathbb{P}_{\mu }\left( Z_{t_{1}}^{\lambda ,\theta }\in
E_{1},...,Z_{t_{n}}^{\lambda ,\theta }\in E_{n}\right) \\
&=&\left( 2\left( \cosh (\lambda \left( \theta +1\right) T)-1\right) \right)
^{\frac{N}{2}}\int_{E_{1}}d\mathsf{x}_{1}...\int_{E_{n}}d\mathsf{x}_{n} \\
&&\times \dprod\limits_{k=2}^{n}g_{\lambda }\left( \mathsf{x}%
_{k},t_{k}-t_{k-1},\mathsf{x}_{k-1}\right) \times \int_{\mathbb{R}^{N}}d%
\mathsf{x}g_{\lambda }\left( \mathsf{x}_{1},t_{1},\mathsf{x}\right)
g_{\lambda }\left( \mathsf{x},\left( \theta +1\right) T-t_{n},\mathsf{x}%
_{n}\right) \\
&=&\left( 2\left( \cosh (\lambda \left( \theta +1\right) T)-1\right) \right)
^{\frac{N}{2}}\int_{E_{1}}d\mathsf{x}_{1}...\int_{E_{n}}d\mathsf{x}_{n} \\
&&\times \dprod\limits_{k=2}^{n}g_{\lambda }\left( \mathsf{x}%
_{k},t_{k}-t_{k-1},\mathsf{x}_{k-1}\right) \times g_{\lambda }\left( \mathsf{%
x}_{1},\left( \theta +1\right) T-(t_{n}-t_{1}),\mathsf{x}_{n}\right)
\end{eqnarray*}%
where we have once again used the semi-group composition law (\ref{semigroup}%
) twice, so that the Gaussian density in $\mathbb{R}^{nN}$ reads%
\begin{eqnarray*}
&&\left( 2\left( \cosh (\lambda \left( \theta +1\right) T)-1\right) \right)
^{\frac{N}{2}} \\
&&\times \dprod\limits_{k=2}^{n}g_{\lambda }\left( \mathsf{x}%
_{k},t_{k}-t_{k-1},\mathsf{x}_{k-1}\right) \times g_{\lambda }\left( \mathsf{%
x}_{1},\left( \theta +1\right) T-(t_{n}-t_{1}),\mathsf{x}_{n}\right) \\
&=&\left( 2\pi \right) ^{-\frac{nN}{2}}\left( \frac{2\lambda ^{n}\left(
\cosh (\lambda \left( \theta +1\right) T)-1\right) }{\sinh \left( \lambda
(\left( \theta +1\right) T-(t_{n}-t_{1})\right) )}\right) ^{\frac{N}{2}%
}\dprod\limits_{k=2}^{n}\left( \sinh (\lambda (t_{k}-t_{k-1})\right) ^{-%
\frac{N}{2}} \\
&&\times \exp \left[ -\frac{\lambda }{2}\sum_{k=2}^{n}\frac{\cosh (\lambda
(t_{k}-t_{k-1}))\left( \left\vert \mathsf{x}_{k}\right\vert ^{2}+\left\vert 
\mathsf{x}_{k-1}\right\vert ^{2}\right) -2\left( \mathsf{x}_{k}\mathsf{,x}%
_{k-1}\right) _{\mathbb{R}^{N}}}{\sinh \left( \lambda (t_{k}-t_{k-1})\right) 
}\right] \\
&&\times \exp \left[ -\frac{\lambda }{2}\frac{\cosh (\lambda (\left( \theta
+1\right) T-(t_{n}-t_{1})))\left( \left\vert \mathsf{x}_{1}\right\vert
^{2}+\left\vert \mathsf{x}_{n}\right\vert ^{2}\right) -2\left( \mathsf{x}_{1}%
\mathsf{,x}_{n}\right) _{\mathbb{R}^{N}}}{\sinh \left( \lambda (\left(
\theta +1\right) T-(t_{n}-t_{1}))\right) }\right] .
\end{eqnarray*}%
For the sake of clarity we now identify the corresponding matrix $C_{\lambda
,\theta }^{-1}$ by considering the case $n=2$ separately from the case $%
n\geq 3$, as $C_{\lambda ,\theta }^{-1}$ is no longer tridagonal. For $n=2$
we obtain%
\begin{equation*}
C_{\lambda ,\theta ,k,k}^{-1}=\frac{\lambda \sinh (\lambda \left( \theta
+1\right) T)}{\sinh \left( \lambda (t_{2}-t_{1}\right) )\sinh \left( \lambda
(\left( \theta +1\right) T-(t_{2}-t_{1})\right) )}\text{ \ \ \ \ for }k=1,2
\end{equation*}%
\bigskip and%
\begin{equation*}
C_{\lambda ,\theta ,2,1}^{-1}=C_{\lambda ,\theta ,1,2}^{-1}=-\frac{\lambda }{%
\sinh \left( \lambda \left\vert t_{2}-t_{1}\right\vert \right) }\text{ }-%
\frac{\lambda }{\sinh \left( \lambda (\left( \theta +1\right) T-\left\vert
t_{2}-t_{1}\right\vert \right) )},
\end{equation*}%
while for $n\geq 3$ we get%
\begin{equation*}
C_{\lambda ,\theta ,k,k}^{-1}=\left\{ 
\begin{array}{c}
\frac{\lambda \sinh \left( \lambda (\left( \theta +1\right)
T-(t_{n}-t_{2})\right) )}{\sinh \left( \lambda (t_{2}-t_{1}\right) )\sinh
\left( \lambda (\left( \theta +1\right) T-(t_{n}-t_{1})\right) )}\text{ \ \
\ \ for }k=1, \\ 
\\ 
\frac{\lambda \sinh \left( \lambda (t_{k+1}-t_{k-1})\right) }{\sinh \left(
\lambda (t_{k+1}-t_{k})\right) \sinh \left( \lambda (t_{k}-t_{k-1})\right) }%
\text{ \ \ \ \ for }k=2,...,n-1, \\ 
\\ 
\frac{\lambda \sinh \left( \lambda (\left( \theta +1\right)
T-(t_{n-1}-t_{1})\right) )}{\sinh \left( \lambda (t_{n}-t_{n-1})\right)
\sinh \left( \lambda (\left( \theta +1\right) T-(t_{n}-t_{1})\right) )}\text{
\ \ \ \ for }k=n\text{ ,}%
\end{array}%
\right.
\end{equation*}%
\begin{equation*}
C_{\lambda ,\theta ,k,k-1}^{-1}=C_{\lambda ,\theta ,k-1,k}^{-1}=-\frac{%
\lambda }{\sinh \left( \lambda \left\vert t_{k}-t_{k-1}\right\vert \right) }%
\text{ \ \ for }k=2,...,n,
\end{equation*}%
and 
\begin{equation*}
C_{\lambda ,\theta ,n,1}^{-1}=C_{\lambda ,\theta ,1,n}^{-1}=-\frac{\lambda }{%
\sinh \left( \lambda (\left( \theta +1\right) T-\left\vert
t_{n}-t_{1}\right\vert \right) )},
\end{equation*}%
all the remaining matrix elements being zero. In both cases we then obtain
by inversion%
\begin{equation*}
C_{\lambda ,\theta ,k,l}=\frac{\sinh \left( \lambda \left\vert
t_{k}-t_{l}\right\vert \right) -\sinh \left( \lambda (\left\vert
t_{k}-t_{l}\right\vert -\left( \theta +1\right) T)\right) }{2\lambda \left(
\cosh (\lambda \left( \theta +1\right) T)-1\right) }
\end{equation*}%
for all $k,l\in \left\{ 1,...,n\right\} $ or, equivalently,%
\begin{equation*}
C_{\lambda ,\theta ,k,l}=\frac{\cosh \left( \lambda \left( \left\vert
t_{k}-t_{l}\right\vert -\frac{\left( \theta +1\right) T}{2}\right) \right) }{%
2\lambda \sinh \left( \frac{\lambda \left( \theta +1\right) T}{2}\right) },
\end{equation*}%
so that (\ref{covariancebis}) eventually holds. Let us now consider the case 
$\theta =0$, namely, the case corresponding to the measure $\hat{\mu}%
_{\lambda ,+}$ defined from (\ref{limitingcase}), for which (\ref%
{variancebis}) and (\ref{covariancebis}) become%
\begin{equation}
\sigma _{\lambda }=\frac{\sinh \left( \lambda T\right) }{2\lambda \left(
\cosh (\lambda T)-1\right) }  \label{variancequarto}
\end{equation}%
and%
\begin{equation}
\mathbb{E}_{\mu }\left( Z_{s}^{\lambda ,+,i}Z_{t}^{\lambda ,+,j}\right) =%
\frac{\cosh \left( \lambda \left( \left\vert t-s\right\vert -\frac{T}{2}%
\right) \right) }{2\lambda \sinh \left( \frac{\lambda T}{2}\right) }\delta
_{i,j}  \label{covariancequarto}
\end{equation}%
for all $i,j\in \left\{ 1,...,N\right\} $, respectively. Let us also
consider the forward It\^{o} integral equation with periodic boundary
conditions%
\begin{eqnarray}
X_{t} &=&e^{-\lambda t}X_{0}+\int_{0}^{t}e^{-\lambda (t-\tau )}dW_{\tau },%
\text{ \ \ }t\in \left[ 0,T\right] ,  \notag \\
X_{0} &=&X_{T}  \label{itoperiodic}
\end{eqnarray}%
rather than (\ref{itoequation}). It is known from Theorem 2.1 in \cite%
{kwakernaak} (see also \cite{pedersen} or Section 5 in \cite{roellythieullen}
for the case $N=1$) that the solution to (\ref{itoperiodic}) can be written
out explicitly and defines a non-Markovian centered Gaussian stationary
process, namely, the so-called periodic Ornstein-Uhlenbeck process given by%
\begin{equation*}
X_{t}=\frac{e^{-\lambda t}}{1-e^{-\lambda T}}\int_{0}^{T}e^{-\lambda (T-\tau
)}dW_{\tau }+\int_{0}^{t}e^{-\lambda (t-\tau )}dW_{\tau },\text{ \ \ }t\in %
\left[ 0,T\right] ,
\end{equation*}%
whose variance and covariance are given by (\ref{variancequarto}) and (\ref%
{covariancequarto}), respectively. Therefore, $Z_{\tau \in \left[ 0,T\right]
}^{\lambda ,+}$ identifies in law with that process. \ \ $\blacksquare $\ \ 

\bigskip

\textsc{Remarks.} (1) It is of course also \textit{a posteriori} clear that
the processes of Proposition 4 are not Markovian since the covariances (\ref%
{covariancebis}) do not factorize as the product of a function of $s$ times
a function of $t$, in contrast to all the cases investigated in Section 3.
Furthermore, a very different way of understanding such factorization
properties in the Markovian case was put forward in Section 6 of \cite%
{kolsrudzambrini}, where the covariances were written as the product of two
linearly independent solutions to some suitable Sturm-Liouville problems.

(2) Problem (\ref{itoperiodic}) is part of a more general class of linear
stochastic differential equations that were investigated by several authors,
including \cite{kwakernaak}, \cite{norris}-\cite{roellythieullen} and some
of the references therein. In this context we ought to mention an analysis
of the law of the solution to (\ref{itoperiodic}) in one dimension carried
out in Section 5 of \cite{roellythieullen}, which establishes a relation
with Bernstein processes whose state space is one-dimensional. The main tool
used there is a formula of integration by parts proved directly on the
underlying infinite-dimensional path space by means of Malliavin's calculus.
This is in sharp contrast to the method we have developed in this section,
as we have first constructed a one-parameter family of non-Markovian, $N$%
-dimensional Bernstein processes associated with the infinite system (\ref%
{forwardharmonicter})-(\ref{backwardharmonicter}), which we have only\textit{%
\ a posteriori }identified with the solution to (\ref{itoperiodic}) when $%
\theta =0$.

(3) While the processes $Z_{\tau \in \left[ 0,T\right] }^{\lambda ,\theta
>0} $ materialize a one-parameter family of non-Markovian random curves in $%
\mathbb{R}^{N}$, they might also occur naturally in completely different
contexts, as the limiting process $Z_{\tau \in \left[ 0,T\right] }^{\lambda
,+}$ does. Thus, this process is quite relevant to the mathematical
investigation of certain quantum systems in equilibrium with a thermal bath
since it identifies with the Gaussian process of mean zero used in Theorem
2.1 of \cite{hoeghkrohn} to compute the expectations of some relevant
physical quantities in statistical mechanics. Indeed, their covariances
coincide since the equality 
\begin{equation}
\frac{\cosh \left( \lambda \left( t-\frac{T}{2}\right) \right) }{2\lambda
\sinh \left( \frac{\lambda T}{2}\right) }=\frac{\exp \left[ -\lambda t\right]
+\exp \left[ -\lambda \left( T-t\right) \right] }{2\lambda \left( 1-\exp %
\left[ -\lambda T\right] \right) }  \label{covequality}
\end{equation}%
holds for every $t\in \left[ 0,T\right] $, the period there being the
inverse temperature, and since the right-hand side of (\ref{covequality}) is
the preferred form to write the covariance in \cite{hoeghkrohn} and \cite%
{kleinlandau}. More generally, certain multidimensional time-periodic
processes such as $Z_{\tau \in \left[ 0,T\right] }^{\lambda ,+}$ have been
useful to solve filtering, smoothing and prediction problems as in \cite%
{kwakernaak}, and also have played an important r\^{o}le in areas as diverse
as the analysis of random evolution of loops as in \cite{norris}, or the
elaboration of global shape models generated by periodic lattices as in
Chapters 11 and 16 of \cite{grenander}, where boundaries of random planar
solid objects were investigated. It is our contention that the whole family
of processes $Z_{\tau \in \left[ 0,T\right] }^{\lambda ,\theta >0}$ might
have an important r\^{o}le to play in those areas as well.

(4) We have already noted that many of the above developments rest on the
series expansion for Green's function which we prove in the appendix below.
We can get expansions of the same kind for Green's functions associated with
forward-backward systems of the form%
\begin{align*}
\partial _{t}u(\mathsf{x},t)& =\frac{1}{2}\Delta _{\mathsf{x}}u(\mathsf{x}%
,t)-\left( \frac{\lambda ^{2}}{2}\left\vert \mathsf{x}\right\vert
^{2}+\kappa W(\mathsf{x})\right) u(\mathsf{x},t),\text{ \ }(\mathsf{x},t)\in 
\mathbb{R}^{N}\times \left( 0,T\right] , \\
u(\mathsf{x},0)& =\varphi _{0}(\mathsf{x}),\text{ \ \ }\mathsf{x}\in \mathbb{%
R}^{N}
\end{align*}%
and%
\begin{align*}
-\partial _{t}v(\mathsf{x},t)& =\frac{1}{2}\Delta _{\mathsf{x}}v(\mathsf{x}%
,t)-\left( \frac{\lambda ^{2}}{2}\left\vert \mathsf{x}\right\vert
^{2}+\kappa W(\mathsf{x})\right) v(\mathsf{x},t),\text{ \ }(\mathsf{x},t)\in 
\mathbb{R}^{N}\times \left[ 0,T\right) , \\
v(\mathsf{x,}T)& =\psi _{T}(\mathsf{x}),\text{ \ \ }\mathsf{x}\in \mathbb{R}%
^{N},
\end{align*}%
for suitable anharmonic potentials $W:\mathbb{R}^{N}\mapsto \mathbb{R}$
where $\kappa >0$. We defer the analysis of Bernstein processes associated
with such systems to a separate publication.

\bigskip

\textbf{Acknowledgments. }The research of the first author was supported in
part by the FCT of the Portuguese government under Grants UID/MAT/04561/2013
and PTDC/MAT/120354/2010, as well as by funds from the group Probabilit\'{e}%
s et Statistiques of the Universit\'{e} de Lorraine. The research of the
second author was also supported by the FCT under Grant
PTDC/MAT/120354/2010. Last but not least, the first author also wishes to
thank the GFMUL and the Complexo Interdisciplinar da Universidade de Lisboa
for their warm hospitality.

\bigskip

We complete this article by proving an important property of (\ref%
{Mehlerkernel}) which we used in the main part of the article.

\section{Appendix: a series expansion for Mehler's N-dimensional kernel}

The notation in this appendix is, of course, the same as in the preceding
sections. The expansion in question is the following:

\bigskip

\textbf{Proposition A1.} \textit{We have }%
\begin{eqnarray}
&&g_{\lambda }(\mathsf{x},t,\mathsf{y})  \notag \\
&=&\sum_{\mathsf{m}\in \mathbb{N}^{N}}e^{-\left( \sum_{j=1}^{N}m_{j}+\frac{N%
}{2}\right) \lambda t}\otimes _{j=1}^{N}h_{m_{j},\lambda }\left( \mathsf{x}%
\right) \times \otimes _{j=1}^{N}h_{m_{j},\lambda }\left( \mathsf{y}\right)
\label{expansion}
\end{eqnarray}%
\textit{where the series converges absolutely for every }$t\in \left( 0,T%
\right] $ \textit{and uniformly for all }$\mathsf{x},\mathsf{y}\in \mathbb{R}%
^{N}$\textit{. Furthermore, (\ref{expansion}) is indeed Green's function
associated with the partial differential equations in (\ref{initial2}).}

\bigskip

\textbf{Proof.} We first prove the result for the one-dimensional case $N=1$%
, namely,%
\begin{eqnarray}
&&\left( \frac{\lambda }{2\pi \sinh \left( \lambda t\right) }\right) ^{\frac{%
1}{2}}\exp \left[ -\frac{\lambda \left( \cosh \left( \lambda t\right) \left(
x^{2}+y^{2}\right) -2xy\right) }{2\sinh \left( \lambda t\right) }\right] 
\notag \\
&=&\sum_{m=0}^{+\infty }e^{-\left( m+\frac{1}{2}\right) \lambda
t}h_{m,\lambda }\left( x\right) h_{m,\lambda }\left( y\right) .
\label{onedimensionalexp}
\end{eqnarray}%
Our point of departure for this is Mehler's notable formula in the form%
\begin{eqnarray*}
&&\left( \frac{1}{\pi \left( 1-\gamma ^{2}\right) }\right) ^{\frac{1}{2}%
}\exp \left[ -\frac{\left( 1+\gamma ^{2}\right) \left( x^{2}+y^{2}\right)
-4\gamma xy}{2\left( 1-\gamma ^{2}\right) }\right] \\
&=&\sum_{m=0}^{+\infty }\gamma ^{m}h_{m}\left( x\right) h_{m}\left( y\right)
\end{eqnarray*}%
valid for every $\gamma \in \left( -1,+1\right) $, which allows one to
express the probability density of two jointly Gaussian variables as a power
series in the correlation parameter $\gamma $ (see, e.g., \cite{slepian}).
By using the scaled Hermite functions (\ref{scaledhermitefunctions})
instead, we obtain%
\begin{eqnarray}
&&\left( \frac{\lambda }{\pi \left( 1-\gamma ^{2}\right) }\right) ^{\frac{1}{%
2}}\exp \left[ -\frac{\lambda \left( \left( 1+\gamma ^{2}\right) \left(
x^{2}+y^{2}\right) -4\gamma xy\right) }{2\left( 1-\gamma ^{2}\right) }\right]
\notag \\
&=&\sum_{m=0}^{+\infty }\gamma ^{m}h_{m,\lambda }\left( x\right)
h_{m,\lambda }\left( y\right) ,  \label{scaledmehler}
\end{eqnarray}%
so that in order to prove (\ref{onedimensionalexp}) we only need to identify 
$\gamma $. To this end we first compare the argument of the exponential on
the left-hand side of (\ref{onedimensionalexp}) with that of (\ref%
{scaledmehler}), which gives the two conditions%
\begin{equation}
\frac{1+\gamma ^{2}}{1-\gamma ^{2}}=\coth \left( \lambda t\right)
\label{firstcondition}
\end{equation}%
and%
\begin{equation}
\frac{\gamma }{1-\gamma ^{2}}=\frac{1}{2\sinh \left( \lambda t\right) }.
\label{secondcondition}
\end{equation}%
Dividing (\ref{firstcondition}) by (\ref{secondcondition}) then leads to the
quadratic equation%
\begin{equation*}
\gamma ^{2}-2\cosh \left( \lambda t\right) \gamma +1=0,
\end{equation*}%
whose only solution in the interval $\left( -1,+1\right) $ is%
\begin{equation}
\gamma =e^{-\lambda t}.  \label{root}
\end{equation}%
The substitution of (\ref{root}) into (\ref{scaledmehler}) then gives (\ref%
{onedimensionalexp}) after some simple algebraic manipulations. Now, using (%
\ref{onedimensionalexp}) we obtain%
\begin{eqnarray*}
&&\sum_{\mathsf{m}\in \mathbb{N}^{N}}e^{-\left( \sum_{j=1}^{N}m_{j}+\frac{N}{%
2}\right) \lambda t}\otimes _{j=1}^{N}h_{m_{j},\lambda }\left( \mathsf{x}%
\right) \times \otimes _{j=1}^{N}h_{m_{j},\lambda }\left( \mathsf{y}\right)
\\
&=&\dprod\limits_{j=1}^{N}\sum_{m_{j}=0}^{+\infty }e^{-\left( m_{j}+\frac{1}{%
2}\right) \lambda t}h_{m_{j},\lambda }\left( x_{j}\right) h_{m_{j},\lambda
}\left( y_{j}\right) \\
&=&\left( \frac{\lambda }{2\pi \sinh \left( \lambda t\right) }\right) ^{%
\frac{N}{2}}\dprod\limits_{j=1}^{N}\exp \left[ -\frac{\lambda \left( \cosh
\left( \lambda t\right) \left( x_{j}^{2}+y_{j}^{2}\right)
-2x_{j}y_{j}\right) }{2\sinh \left( \lambda t\right) }\right]
\end{eqnarray*}%
which is (\ref{expansion}), and the uniform convergence of the series
follows from the Cr\'{a}mer-Charlier inequality which ensures the uniform
boundedness of the $h_{m}$'s in $m$ and in their argument (see, e.g.,
Section 10.18 in \cite{erdmagobertri} and the references therein). Finally,
using the series expansion (\ref{expansion}) and the eigenvalue equation (%
\ref{eigenequation}) of Lemma 1 we get%
\begin{eqnarray*}
\partial _{t}g_{\lambda }(\mathsf{x},t,\mathsf{y}) &=&\frac{1}{2}\Delta _{%
\mathsf{x}}g_{\lambda }(\mathsf{x},t,\mathsf{y})-\frac{\lambda
^{2}\left\vert \mathsf{x}\right\vert ^{2}}{2}g_{\lambda }(\mathsf{x},t,%
\mathsf{y}), \\
(\mathsf{x},t) &\in &\mathbb{R}^{N}\mathbb{\times }\left( 0,T\right]
\end{eqnarray*}%
for every $\mathsf{y}\in \mathbb{R}^{N}$, and%
\begin{equation*}
\lim_{t\rightarrow 0_{+}}g_{\lambda }(\mathsf{x},t,\mathsf{y})=\sum_{\mathsf{%
m}\in \mathbb{N}^{N}}\otimes _{j=1}^{N}h_{m_{j},\lambda }\left( \mathsf{x}%
\right) \times \otimes _{j=1}^{N}h_{m_{j},\lambda }\left( \mathsf{y}\right)
=\delta (\mathsf{x-y)}
\end{equation*}%
in the sense of distributions since the $\otimes _{j=1}^{N}h_{m_{j},\lambda
} $'s constitute a complete orthonormal system in $L^{2}\left( \mathbb{R}%
^{N},\mathbb{C}\right) $. Consequently, $g_{\lambda }$ is indeed Green's
function associated with the partial differential equation in (\ref{initial2}%
). $\ \ \blacksquare $

\bigskip

\textsc{Remark.} Of course, the semi-group composition law%
\begin{eqnarray}
&&g_{\lambda }(\mathsf{x},s+t,\mathsf{y})  \notag \\
&=&\int_{\mathbb{R}^{N}}d\mathsf{z}g_{\lambda }(\mathsf{x},s,\mathsf{z}%
)g_{\lambda }(\mathsf{z},t,\mathsf{y})  \label{semigroup}
\end{eqnarray}%
valid for all $\mathsf{x},\mathsf{y}\in \mathbb{R}^{N}$ is inherent in the
fact that $g_{\lambda }$ is Green's function for (\ref{initial2}), but
follows most directly from the series expansion (\ref{expansion}) and the
orthogonality relations%
\begin{equation*}
\int_{\mathbb{R}^{N}}d\mathsf{x}\otimes _{j=1}^{N}h_{m_{j},\lambda }\left( 
\mathsf{x}\right) \times \otimes _{j=1}^{N}h_{n_{j},\lambda }\left( \mathsf{x%
}\right) =\dprod\limits_{j=1}^{N}\delta _{m_{j},n_{j}}.
\end{equation*}%
Furthermore, we also observe that (\ref{expansion}) brings out the entire
pure point spectrum of Lemma 1 in the argument of the exponential, a fact
that is crucial in our construction of Section 4.

\end{document}